\documentclass{article}
\pdfoutput=1

\usepackage[nonatbib]{arxiv}
\usepackage[utf8]{inputenc} 
\usepackage[T1]{fontenc}    
\usepackage{hyperref}       
\usepackage{url}            
\usepackage{booktabs}       
\usepackage{amsfonts}       
\usepackage{nicefrac}       

\usepackage[title]{appendix}
\usepackage{IEEEtrantools}
\usepackage{amsmath}
\usepackage{amssymb}
\usepackage{amsthm}
\usepackage{mathtools}
\usepackage{siunitx}
\usepackage[hang,flushmargin]{footmisc} 

\newtheorem{theorem}{Theorem}
\newtheorem{lem}{Lemma}
\newtheorem{assum}{Assumption}

\newtheorem{corol}{Corollary}
\newcommand{\cmag}[1]{\vert #1 \vert}
\newcommand{\thedoi}{10.1109/TCST.2020.3044862}

\let\svthefootnote\thefootnote
\newcommand\blankfootnote[1]{%
  \let\thefootnote\relax\footnotetext{#1}%
  \let\thefootnote\svthefootnote%
}

\begin{document}
\title{Model-based resonance tracking of linear systems}
\rhead{\scshape DOI: \href{www.doi.org/\thedoi}{\thedoi} - \today}
\author{Thomas Vasileiou}
\date{\today}
\blankfootnote{
\noindent This  paper  has  been  accepted  for  publication  in IEEE Transactions on Control Systems Technology.\\
DOI: \href{www.doi.org/\thedoi}{\thedoi}  \\
\textsuperscript{\copyright}~2020 IEEE.  Personal use of this material is permitted.  Permission from IEEE must be obtained for all other uses, in any current or future media, including reprinting/republishing this material for advertising or promotional purposes, creating new collective works, for resale or redistribution to servers or lists, or reuse of any copyrighted component of this work in other works.
}
\maketitle

\begin{abstract}
The present paper develops recursive algorithms to track shifts in the resonance frequency of linear systems in real time.
To date, automatic resonance tracking has been limited to non-model-based approaches, which rely solely on the phase difference between a specific input and output of the system.
Instead, we propose a transformation of the system into a complex-valued representation, which allows us to abstract the resonance shifts as an exogenous disturbance acting on the excitation frequency, perturbing the excitation frequency from the natural frequency of the plant.
We then discuss the resonance tracking task in two parts: recursively identifying the frequency disturbance and incorporating an update of the excitation frequency in the algorithm.
The complex representation of the system simplifies the design of resonance tracking algorithms due to the applicability of well-established techniques.
We discuss the stability of the proposed scheme, even in cases that seriously challenge current phase-based approaches, such as nonmonotonic phase differences and multiple-input multiple-output systems.
Numerical simulations further demonstrate the performance of the proposed resonance tracking scheme.

\end{abstract}


\section{Introduction}
Precisely tracking the resonance frequency of oscillating systems is of great interest in resonant sensing \cite{Hauptmann1991,Boisen2011} and in the driving of vibrating loads \cite{Claeyssen2003,Gokcek2005}.
Resonant sensors, the function of which relies on the resonant characteristic of a vibrating structure, have been proposed for a wide range of measurements and instruments, including thermometers \cite{Larsen2011}, accelerometers \cite{Aikele2001}, viscometers \cite{Brack2016}, humidity sensors \cite{Sheng2011}, water cut measurements \cite{Avila2017} and gyroscopes \cite{Raman2009,Fei2015}.
In terms of miniaturization and increased sensitivity, microelectromechanical systems (MEMS) with vibrating cantilevers have emerged as an appealing solution, and achievements such as atomic force microscopy in space \cite{Bentley2016} and mass detection in the range of atto- and zeptograms \cite{Boisen2011} have been reported.
Furthermore, resonant electromechanical actuators have been widely proposed for power electronics \cite{Li2016,Singh2014,Bosshard2014,Choi2006}, ultrasonic applications \cite{Liu2015}, thermosonic wire bonding \cite{Zhang2015} and acoustic particle trapping \cite{Hammarstroem2014}.

To increase the sensitivity of sensors and the power output of vibrating actuators, designers adopt systems with ``sharp'' resonances (low damping and a high quality factor) \cite{Hammarstroem2014}. As a drawback, this design leads to diminished performance when the excitation frequency deviates even slightly from the resonance frequency due to the inherently narrow bandwidth of the system.
Even for actuators that have been designed to operate at a constant resonance, shifts from the designed operating frequency may occur because of environmental changes such as temperature and humidity variations \cite{Ferguson2005}, aging of the device \cite{Sun2002} or changes in the load \cite{Zhang2015,Bosshard2014}. As a remedy, designers resort to feedback resonance tracking control to compensate for these shifts and achieve maximum efficiency \cite{Gokcek2005,Li2016}.
In the case of sensing applications based on changes in the resonance frequency with the measured quantity, the use of feedback control is unavoidable, and the performance of the control system directly affects the sensitivity, resolution and bandwidth of the sensor \cite{Sun2002,Albrecht1991}.

Regardless of the application, the phase locked loop (PLL) is the typical scheme for resonance tracking \cite{Sun2002,Brack2016}.
The self-sustained oscillation (SSO) scheme \cite{Lulec2016,Gorreta2016} has also been proposed to induce excitation at resonance frequencies.
Both techniques achieve tracking of the resonance frequency by maintaining a constant phase difference between the input and the output of the system. Their main difference lies in the fact that the PLL utilizes an external oscillator to generate the signal that excites the system. In contrast, for the SSO case, the excitation signal is generated by the oscillating structure itself; the system output is amplified and phase-shifted before it is fed back to the system. The design and analysis of the PLL and SSO resonant tracking schemes connected with a single-input single-output (SISO) 2\textsuperscript{nd}-order system have been extensively discussed by many authors for various applications \cite{Sun2002,Park2009,Aikele2001,SohanianHaghighi2012,Lulec2016,Brack2016}.
A drawback of these resonance tracking approaches is that they are far from being model-based; more importantly, however, the closed-loop robustness and stability cannot be guaranteed.
This situation is true even for linear SISO systems if the phase difference between the input and output signal is nonmonotonic. This problem has been pointed out for piezoelectric actuators and multi-degree-of-freedom systems, where resonance and anti-resonance frequencies are present \cite{Hayashi1992,Zhang2015,Brack2016a,Liu2015}.

To supplement the two techniques mentioned above, specialized resonance tracking algorithms have been developed for specific cases. 
In the absence of phase information, an algorithm that detects the maximum of the output signal has been employed in piezoresponse force microscopy \cite{Rodriguez2007}.
A maximum power tracking adaptive approach for the driving of resonant loads was proposed in \cite{Gokcek2005}; in this approach, a small sinusoidal perturbation signal is added to the driving signal to estimate the derivative of the absorbed power and update the excitation frequency.
A control algorithm that tunes the system to a specific resonance frequency was developed to address the issue of online modal frequency matching in vibratory gyroscopes \cite{Park2009}.
Other specialized schemes employ controller scheduling \cite{Song2016} and fuzzy logic \cite{Zhang2015}.
Nonetheless, all of the previous approaches were developed for operation with a specific system, which makes it difficult to generalize the approach and may require multiple driving signals.
Moreover, the application of these schemes has been limited to SISO systems, making their extension to multiple-input multiple-output (MIMO) systems rather complicated.

The PLL and SSO techniques have both proven their capabilities in numerous applications. By neglecting the system model, these techniques are sufficiently general and therefore applicable in many cases.
On the other hand, tuning the controller parameters may be tedious, and theoretical analysis of the closed-loop performance still requires a mathematical description of the system.
Nevertheless, if a description of the system is available, a model-based controller design approach is expected to yield improved performance.
Here, we consider control algorithms that are able to track the resonance of linear models.
We introduce a representation of the oscillating plant, which splits the system into a linear time-invariant (LTI) model and an abstract shift of the frequency exciting the system.
This model representation simplifies the application of well-established control and estimation techniques.

The remainder of this paper is organized as follows. The resonance tracking problem is formalized in Section~\ref{sec:problem}. The complex representation of oscillating linear systems and the abstract modeling of resonance shifts are introduced in Section~\ref{sec:complexss}. Section~\ref{sec:freqtracking} presents the estimation of the resonance shift and the update of the excitation frequency. Implementation considerations and numerical simulations are discussed in Section~\ref{sec:simulations}. Section~\ref{sec:conclusion} concludes the paper.

\section{Problem statement}\label{sec:problem}
We consider the following discrete-time linear model:
\begin{IEEEeqnarray}{rCl}
x_{k+1} & = & \tilde{A} x_k + \tilde{B} u_k + \tilde{w}_k \label{eq:linear_uncertain_state} \\
y_k & = & \tilde{C} x_k + \tilde{D} u_k + \tilde{v}_k \label{eq:linear_uncertain_output}
\end{IEEEeqnarray}
where $x_k \in \mathbb{R}^n$, $u_k \in \mathbb{R}^m$, and $y_k \in \mathbb{R}^p$ denote the state, input and output vectors, respectively, at the discrete sampling instances $k \in \mathbb{N}$ and $\tilde{w}_k  \sim  \mathcal{N}(0,Q)$ and $\tilde{v}_k \sim \mathcal{N}(0,R)$ are uncorrelated additive white Gaussian noise used to model the disturbance input and the measurement noise, respectively.
The system matrices depend on an unknown parameter vector, $\kappa$, such as $\tilde{A}(\kappa)$, $\tilde{B}(\kappa)$, $\tilde{C}(\kappa)$ and $\tilde{D}(\kappa)$.
The vector $\kappa$ may be time-varying: in such cases, the system (\ref{eq:linear_uncertain_state}-\ref{eq:linear_uncertain_output}) becomes linear time-varying (LTV), which we indicate by explicitly adding the subscript $k$ to the system matrices, e.g., $\tilde{A}_k$.

\begin{figure}
\centering
\includegraphics{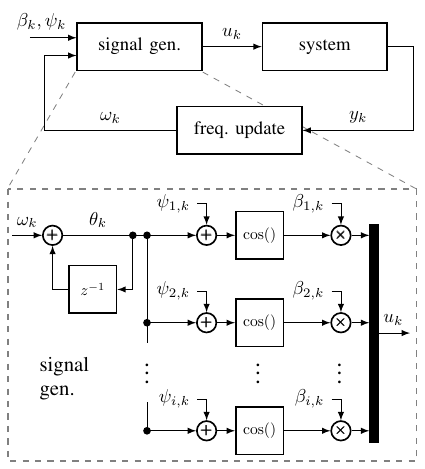}
\caption{Schematic of the frequency tracking problem.}
\label{fig:blockdiagram}
\end{figure}

We assume that the system has at least one pair of conjugated complex eigenvalues corresponding to the resonance frequency of interest.
We denote the eigenvalue of interest by $\tilde{\lambda}_s$ and its corresponding resonance by $\tilde{\omega}_\lambda = \arg\tilde{\lambda}_s > 0$.
The system is subjected to a sinusoidal input $u_{i,k} = \beta_{i,k} \cos(\theta_k + \psi_{i,k})$, where $\beta_{i,k} \in \mathbb{R}_{\ge 0}$ and $\psi_{i,k} \in (-\pi, \, \pi]$ denote the instantaneous amplitude and phase, respectively, of $u_{i,k}$, the $i$\textsuperscript{th} component of $u_k$.
The common reference phase is denoted by $\theta_k \in \mathbb{R}$ and evolves according to
\begin{equation}\label{eq:theta}
\theta_{k+1} = \theta_k + \omega_k
\end{equation}
where $\omega_k$ is the normalized angular frequency of the excitation.

The goal of this study is to develop a recursive scheme that detects the resonance frequency of the linear system~(\ref{eq:linear_uncertain_state}-\ref{eq:linear_uncertain_output}).
More precisely, we aim to develop a recursive algorithm that drives the excitation frequency $\omega_k$ towards $\tilde{\omega}_\lambda$.
A schematic representation of the structure of the resonance tracking problem and the composition of $u_k$ is shown in Fig.~\ref{fig:blockdiagram}.

\section{Complex state-space model}\label{sec:complexss}
In this section, we introduce a transformation of the oscillating system into an equivalent description where the state, input and output variables are represented as complex envelopes of sinusoidal signals.

\subsection{CSS representation}
We discuss the transformation for a general LTV system:
\begin{IEEEeqnarray}{rCl}
x_{k+1} & = & \tilde{A}_k x_k + \tilde{B}_k u_k + \tilde{w}_k \label{eq:ltv_state} \\
y_k & = & \tilde{C}_k x_k + \tilde{D}_k u_k + \tilde{v}_k \label{eq:ltv_output}
\end{IEEEeqnarray}
excited by the sinusoidal input $u_k$.
Inspired by \cite{Mcloskey1999, Sun2002}, we write the state and output variables in an amplitude-phase formulation, $x_{i,k} = \alpha_{i,k} \cos ( \theta_k + \phi_{i,k})$ and $y_{i,k} = \delta_{i,k} \cos (\theta_k + \varphi_{i,k})$, where $\alpha_{i,k}$, $\delta_{i,k}$, $\phi_{i,k}$ and $\varphi_{i,k}$ are components of time-varying vectors of appropriate dimensions.
Substituting the phase-amplitude expressions into~(\ref{eq:ltv_state}), the time update for the $i$\textsuperscript{th} component of $x_k$ is
\begin{equation}
\alpha_{i,k+1} \cos (\theta_k + \omega_k + \phi_{i,k+1}) =  \sum_{l=1}^n \tilde{A}_{il,k} \alpha_{l,k} \cos (\theta_k + \phi_{l,k}) \IEEEnonumber + \sum_{l=1}^m \tilde{B}_{il,k} \beta_{l,k} \cos (\theta_k + \psi_{l,k}) \label{eq:cos_state}
\end{equation}
where $\tilde{A}_{il,k}$ and $\tilde{B}_{il,k}$ are the components of the $i$\textsuperscript{th} row and $l$\textsuperscript{th} column of the matrices $\tilde{A}_k$ and $\tilde{B}_k$, respectively.
In the previous expressions, we have neglected the effect of $\tilde{w}_k$.
We use the angle-sum trigonometric identities to expand the terms in~(\ref{eq:cos_state}):
\begin{IEEEeqnarray}{rCl}
\IEEEeqnarraymulticol{3}{l}{\alpha_{i,k+1} [ \cos \theta_{k} \cos(\omega_k + \phi_{i,k+1}) - \sin \theta_{k} \sin(\omega_k + \phi_{i,k+1}) ] } \IEEEnonumber \\
\quad & = & \sum_{l=1}^n \tilde{A}_{il,k} [ \alpha_{i,k} ( \cos \theta_k \cos \phi_{l,k} - \sin \theta_k \sin \phi_{l,k} ) ] + \sum_{l=1}^m \tilde{B}_{il,k} [ \beta_{l,k} ( \cos \theta_k \cos \psi_{l,k} - \sin \theta_k \sin \psi_{l,k} ) ] \IEEEeqnarraynumspace \label{eq:cos_state_expanded}
\end{IEEEeqnarray}
which can be compactly written as
\begin{equation}
c_{c,i,k} \cos \theta_k + c_{s,i,k} \sin \theta_k =
\sqrt{c_{c,i,k}^2 + c_{s,i,k}^2} \sin (\theta_k + \arctan\frac{c_{c,i,k}}{c_{s,i,k}}) = 0. \label{eq:cos_state_compact}
\end{equation}
Eq.~(\ref{eq:cos_state_compact}) becomes independent of $\theta_k$ by setting $c_{c,i,k} = c_{s,i,k} = 0$:
\begin{IEEEeqnarray}{rCl}
\alpha_{i,k+1} \cos(\omega_k + \phi_{i,k+1}) & = & \sum_{l=1}^n \tilde{A}_{il,k} \alpha_{l,k} \cos \phi_{l,k} + \sum_{l=1}^m \tilde{B}_{il,k} \beta_{l,k} \cos \psi_{l,k} \label{eq:system_cos} \\
\alpha_{i,k+1} \sin(\omega_k + \phi_{i,k+1}) & = & \sum_{l=1}^n \tilde{A}_{il,k} \alpha_{l,k} \sin \phi_{l,k} + \sum_{l=1}^n \tilde{B}_{il,k} \beta_{l,k} \sin \psi_{l,k} \label{eq:system_sin}.
\end{IEEEeqnarray}
We introduce the transformation for the $i$\textsuperscript{th} state component in complex notation, $z_{i,k} = \alpha_{i,k} ( \cos \phi_{i,k} + j \sin \phi_{i,k} ) = \alpha_{i,k} e^{j \phi_{i,k}}$.
Similarly, we write the output as $q_{i,k} = \delta_{i,k} e^{j \varphi_{i,k}}$ and the input as $s_{i,k} = \beta_{i,k} e^{j \psi_{i,k}}$.
The symbol $j = \sqrt {-1}$ is the imaginary unit.
This complex signal form is similar to the complex envelope representation of a bandpass signal in communication channels \cite{Neeser1993} and to the analytic signal \cite{Boashash1992}.
Substitution of the complex variables into~(\ref{eq:system_cos}-\ref{eq:system_sin}) results in the following complex system:
\begin{IEEEeqnarray}{rCl}
z_{k+1} & = & ( \tilde{A}_k z_k + \tilde{B}_k s_k + w_k ) e^{-j \omega_k} \label{eq:linear_complex_state} \\
q_k & = & \tilde{C}_k z_k + \tilde{D}_k s_k + v_k \label{eq:linear_complex_output}
\end{IEEEeqnarray}
where $w_k$ and $v_k$ are proper random variables with a complex Gaussian ($\mathcal{CN}$) probability density function, that is, $w_k \sim \mathcal{CN}(0, \, Q)$ and $v_k \sim \mathcal{CN}(0, \, R)$.
The complex envelope of white real-valued Gaussian signals has been shown to be complex proper normal, where the properness arises from the stationarity assumption \cite{Picinbono1994,Neeser1993}.
Eq.~(\ref{eq:linear_complex_output}) follows from~(\ref{eq:ltv_output}) using the same procedure.
We refer to the system~(\ref{eq:linear_complex_state}-\ref{eq:linear_complex_output}) as the complex state space (CSS) representation.

The conversion of an LTV system into the CSS representation can be derived by substituting the analytic signal directly into~(\ref{eq:ltv_state}-\ref{eq:ltv_output}).
The derivation is simpler and directly relates each signal to its complex envelope but lacks intuition and reasoning for the complex representation of the variables.
The substitution of $u_k = s_k e^{j \theta_k}$ and $\tilde{w}_k = w_k  e^{j \theta_k}$ into (\ref{eq:ltv_state}) results in
\begin{equation}\label{eq:lti_complex_envelope}
x_{k+1} = \tilde{A}_k x_k + \tilde{B} s_k e^{j \theta_k} + w_k  e^{j \theta_k}
\end{equation}
which is equal to (\ref{eq:linear_complex_state}) for $x_k = z_k e^{j \theta_k}$ and $x_{k+1} = z_{k+1} e^{j \omega_k} e^{j \theta_k}$.
The evolution of the real part of (\ref{eq:lti_complex_envelope}), which is entirely disconnected from the complex part, coincides with (\ref{eq:ltv_state}).
Furthermore, we apply the same reasoning to derive the continuous-time CSS representation in Appendix~\ref{app:continuousTimeCSS}.


\subsection{Properties of the CSS representation}
Transforming (\ref{eq:ltv_state}-\ref{eq:ltv_output}) into the CSS representation constitutes an alternative description of the original system.
Under equivalent input and noise sequences, the trajectories of $x_k$ and $z_k$ can be derived from one another.
It is therefore expected that the properties of~(\ref{eq:ltv_state}-\ref{eq:ltv_output}) are retained after the transformation.
In the following, the relevant properties for the design of the resonance tracking algorithm are discussed.

The state transition matrices of the two systems are closely related. Let $\Phi(k,i) = \prod_{l=i}^{k-1} \tilde{A}_l$ be the state transition matrix for the system~(\ref{eq:ltv_state}-\ref{eq:ltv_output}); then, the state transition matrix for the CSS representation (\ref{eq:linear_complex_state}-\ref{eq:linear_complex_output}) is
\begin{equation}
\Phi_c(k,i) = \prod_{l=i}^{k-1} (\tilde{A}_l e^{-j \omega_l}) =  \left (\prod_{l=i}^{k-1} \tilde{A}_l \right ) e^{-j \sum_{l=i}^{k-1} \omega_l}  = \Phi(k,i) e^{-j \sum_{l=i}^{k-1} \omega_l}\label{eq:state_trans_mat}
\end{equation}
with $\Phi(k,k) = \Phi_c(k,k) = I$.

\begin{lem}\label{lem:stb}
Consider the LTV system~(\ref{eq:ltv_state}-\ref{eq:ltv_output}) and the associated CSS representation (\ref{eq:linear_complex_state}-\ref{eq:linear_complex_output}).
For any given sequence $\omega_i, \, i = 1, 2, \dotsc, k$, the system (\ref{eq:linear_complex_state}-\ref{eq:linear_complex_output}) is globally exponentially stable if and only if the associated system (\ref{eq:ltv_state}-\ref{eq:ltv_output}) is globally exponentially stable.
\end{lem}
\begin{proof}
The exponential stability of (\ref{eq:ltv_state}-\ref{eq:ltv_output}) implies that there exist scalars $c_1 > 0$ and $0 < \lambda_\Phi < 1$ such that \cite[lemma~1]{Zhou2017}
\begin{equation}
\lVert \Phi(k,i) \rVert \leq c_1 \lambda_\Phi^{k-i}, \quad \forall i,k \in \mathbb{N}.
\end{equation}
For the CSS representation, $\lVert \Phi_c(k,i) \rVert = \lVert \Phi(k,i) \rVert \leq c_1 \lambda_\Phi^{k-i} $, which concludes the global exponential stability of~(\ref{eq:linear_complex_state}-\ref{eq:linear_complex_output}).

The reverse statement can be shown similarly.
\end{proof}

Furthermore, the optimal control and optimal estimation problems for (\ref{eq:linear_complex_state}-\ref{eq:linear_complex_output}) and the LTV system are directly connected.
Assume an observable system (\ref{eq:ltv_state}-\ref{eq:ltv_output}), and consider the optimal state observer design problem with an initial state estimate $\hat{x}_0 = E [x_0]$ and variance $P_{x,0} = E[(x_0 - \hat{x}_0)^T  (x_0 - \hat{x}_0)]$, where $E[ \cdot ]$ denotes the expected value.
The trajectory of the optimal state estimates and the covariance matrix, $P_{x,k+1|k}$, are given by the Kalman filter equations.
Specifically, $P_{x,k+1|k}$ follows the Riccati equation
\begin{equation}
P_{x,k+1|k} = \tilde{A}_k P_{x,k|k-1} \tilde{A}_k^T + Q - \tilde{A}_k P_{x,k|k-1} \tilde{C}_k^T (\tilde{C}_k P_{x,k|k-1} \tilde{C}_k^T + R)^{-1} \tilde{C}_k P_{x,k|k-1} \tilde{A}_k^T. \label{eq:kalm_ltv}
\end{equation}

\begin{lem}\label{lem:kalman_filter}
Given the optimal estimation problem for the CSS system (\ref{eq:linear_complex_state}-\ref{eq:linear_complex_output}) with the initial conditions $\hat{z}_0 = E [z_0]$ and $P_{z,0} = E[(z_0 - \hat{z}_0)^H (z_0 - \hat{z}_0)]$, (\ref{eq:kalm_ltv}) then describes the evolution of the optimal covariance, $P_{z,k+1|k}$, for any sequence $\omega_i, \, i = 1, 2, \dotsc, k$.
Moreover, if $P_{z,0} = P_{x,0}$, then $P_{z,k+1|k} = P_{x,k+1|k}, \, \forall k > 0$.
\end{lem}
\begin{proof}
The system matrices of (\ref{eq:linear_complex_state}-\ref{eq:linear_complex_output}) are real-valued, and the noise variables are proper.
The optimal estimator for such systems has been shown to be the Kalman filter, which achieves optimality in terms of being unbiased and having minimum variance \cite[remark 6]{Dini2012}.

We write the Kalman filter equations for the CSS model, presented in a prediction and a correction step, as
\begin{IEEEeqnarray}{rCl}
\IEEEeqnarraymulticol{3}{s}{Prediction step:} \IEEEnonumber \\*
\hat{z}_{k|k-1} & = & ( \tilde{A}_k \hat{z}_{k-1|k-1} + \tilde{B}_k s_k ) e^{-j \omega_k} \\
P_{z,k|k-1} & = & (\tilde{A}_k e^{-j \omega_k}) P_{z,k-1|k-1} (\tilde{A}_k e^{-j \omega_k})^H + Q = \tilde{A}_k P_{z,k-1|k-1} \tilde{A}_k^T + Q \label{eq:kalman_pred_var}  \\
\IEEEeqnarraymulticol{3}{s}{Correction step:} \IEEEnonumber \\*
\hat{z}_{k|k} & = & \hat{z}_{k|k-1} + L_k ( q_k - \tilde{C}_k \hat{z}_{k|k-1} - \tilde{D}_k s_k )  \IEEEeqnarraynumspace \\
L_k & = & P_{z,k|k-1} \tilde{C}_k^T ( \tilde{C}_k P_{z,k|k-1} \tilde{C}_k^T + R )^{-1} \label{eq:kalman_gain} \\
P_{z,k|k} & = & (I - L_k \tilde{C}_k) P_{z,k|k-1} \label{eq:kalman_corr_var}
\end{IEEEeqnarray}
where $\hat{z}_{k|k-1}$ and $\hat{z}_{k|k}$ are the prior and posterior state estimates, respectively.
We use the same convention for $P_{z,k+1|k}$ and $P_{z,k|k}$.
Combining (\ref{eq:kalman_corr_var}-\ref{eq:kalman_gain}) with (\ref{eq:kalman_pred_var}) results in (\ref{eq:kalm_ltv}).
Additionally, if $P_{z,0} = P_{x,0}$, the trajectories of $P_{x,k+1|k}$ and $P_{z,k+1|k}$ are identical and independent of $\omega_k$.
\end{proof}

The previous result is connected to the fact that the properness (with respect to the vanishing pseudocovariance) and normality of the complex random variable are retained under affine transformations \cite{Neeser1993}.
Moreover, for the proper random variable, $z_k$, the random variable $z_k e^{-j \omega_k}$ has the same second-order statistical properties \cite{Picinbono1995}.
We note that $P_{z,k+1|k}$ is bounded, a statement that extends the variance of proper random variables filtered by a CSS model.
The connection between the LTV system and the CSS representation can be extended to the linear quadratic regulator (LQR) problem; for more details, see Appendix~\ref{app:optimalLQR}.


\subsection{Modeling resonance shifts}
The model-based resonance tracking algorithm employs a simplified model to control the plant; the real system is abstracted into a nominal LTI part and a complex-valued frequency disturbance, $h$.
The frequency disturbance acts independently to shift the resonance frequency of the system by $\arg h$.
The disturbed system can be written as
\begin{IEEEeqnarray}{rCl}
z_{k+1} & = & h \left( A z_k + B s_k + w_k \right) e^{-j \omega_k} \label{eq:lti_ch_state} \\
q_k & = & C z_k + D s_k + v_k. \label{eq:lti_ch_output}
\end{IEEEeqnarray}
The description (\ref{eq:lti_ch_state}-\ref{eq:lti_ch_output}) is sufficiently general and can model a set of LTI systems with similar dynamics as an average LTI model disturbed by some $h$ or split an LTV system into an LTI part and a time-varying $h_k$.

For a fixed $h$, the model (\ref{eq:lti_ch_state}-\ref{eq:lti_ch_output}) is connected to an equivalent LTI representation, as in (\ref{eq:lti_complex_envelope}).
In this case, the complex and real parts of the system variables are not disconnected but can be written in real-valued form as
\begin{IEEEeqnarray}{rCl}
\begin{bmatrix} \Re x_{k+1} \\ \Im x_{k+1} \end{bmatrix} & = &
\begin{bmatrix*}[r] h_r A & -h_i A \\ h_i A & h_r A \end{bmatrix*}
\begin{bmatrix} \Re x_k \\ \Im x_k  \end{bmatrix} +
\begin{bmatrix*}[r] h_r B & -h_i B \\ h_i B & h_r B \end{bmatrix*}
\begin{bmatrix} \Re u_k \\ \Im u_k \end{bmatrix} \label{eq:lti_realch_state} \\
\begin{bmatrix} \Re y_k \\ \Im y_k \end{bmatrix} & = &
\begin{bmatrix} C & 0 \\ 0 & C \end{bmatrix}
\begin{bmatrix} \Re x_k \\ \Im x_k  \end{bmatrix} +
\begin{bmatrix} D & 0 \\ 0 & D \end{bmatrix}
\begin{bmatrix} \Re u_k \\ \Im u_k \end{bmatrix} \label{eq:lti_realch_output}
\end{IEEEeqnarray}
where $h_r = \Re h$ and $h_i = \Im h$.
In the following, we adopt (\ref{eq:lti_realch_state}-\ref{eq:lti_realch_output}) whenever it is advantageous to include $\omega_k$ in the excitation signal $u_k$; alternatively, we use the CSS representation (\ref{eq:lti_ch_state}-\ref{eq:lti_ch_output}) when it is more convenient to treat $\omega_k$ as a system parameter.

It is also of interest to define a one-step-ahead predictor for (\ref{eq:lti_ch_state}-\ref{eq:lti_ch_output}).
As a consequence of Lemma~\ref{lem:kalman_filter}, such a predictor can be realized as a steady-state Kalman filter:
\begin{IEEEeqnarray}{rCl}
\hat{z}_{k+1} & = & h [ (A - L C) \hat{z}_k + (B - L D) s_k + L q_k ] e^{-j \omega_k} \IEEEeqnarraynumspace \label{eq:css_est_state} \\
\hat{q}_k & = & C \hat{z}_k + D s_k \label{eq:css_est_output}
\end{IEEEeqnarray}
where $\hat{z}_k$ and $\hat{q}_k$ denote the state and output estimates at time $k$, respectively.
The steady-state Kalman gain is
\begin{equation}\label{eq:lti_kalman_gain}
L_{C} = e^{-j \omega_k} A P_\infty C^T (C P_\infty C^T + R)^{-1} = e^{-j \omega_k} L
\end{equation}
where $P_\infty$ is the fixed point for $P_{z,k+1|k}$ in (\ref{eq:kalm_ltv}).
In the following, we denote the set of parameterized models or the parameterized one-step-ahead predictors with $\mathcal{M}$.

\section{Resonance tracking control}\label{sec:freqtracking}

The resonance tracking algorithm developed in the present section is inspired by the ideas of adaptive control.
Initially, we apply the recursive identification method to estimate $h$ without feedback on $\omega_k$, which provides the base for the closed-loop algorithm.
Then, we modify the algorithm to include the update of the excitation frequency.

\subsection{Estimating the frequency shift}

The estimation of $h$ is based on minimizing the estimation error, namely, the discrepancy between the observed and estimated values, $\rho_k = q_k - \hat{q}_k$.
The model identification of dynamical systems based on parameterized estimators has been thoroughly discussed in the literature \cite{Ljung1978,Ljung1981}; the identification technique is termed the prediction error method (PEM).
We follow \cite{Ljung1981} and discuss the recursive PEM in two steps; initially, we transfer previous results for the offline identification procedure to the current complex-valued structure, and then we discuss the recursive counterpart of the PEM.

In the offline setup, the sequences $q_k$, $s_k$ and $\omega_k$ for $k = 1, 2, \dotsc, N$ are available for the identification procedure.
For any given value of $h$, the sequence $\hat{q}_k$ is computed by
(\ref{eq:css_est_state}-\ref{eq:css_est_output}), and the cumulative estimation error is computed by
\begin{equation}\label{eq:estObjFun}
J_{est,N}(h) = \frac{1}{2 N} \sum_{k = 1}^{N} \lVert \rho_k \rVert^2.
\end{equation}
We define $\hat{h} = \arg\!\min J_{est,N}$.
The objective function $J_{est,N}$ is a real-valued function of complex arguments; optimization problems involving such functions have been addressed by the mathematical framework of the $\mathbb{CR}$-calculus \cite{KreutzDelgado2009,Musella2016}.
Although $J_{est,N}$ is not holomorphic, the second-order expansion exists and can form the basis for the Gauss--Newton gradient descent minimization method. The offline estimate is obtained by iteratively applying
\begin{equation}
    \hat{h}^{(i+1)} = \hat{h}^{(i)} + \delta\hat{h}^{(i)}
\end{equation}
from some initial guess $\hat{h}^{(0)}$, where the superscript $(i)$ denotes the values of the procedure at iteration $i$ and $\delta \hat{h}^{(i)}$ is the search direction.
To apply this theory to our problem, we rewrite~(\ref{eq:estObjFun}) as
\begin{equation}
J_{est,N} = \frac{1}{2 N} \sum_{k = 1}^{N} \lVert q_k - C \hat{z}_k - D s_k \rVert^2 = \frac{1}{2 N} [ q_m -  g (h) ]^H [ q_m -  g (h) ] \label{eq:estObjFun_exp}
\end{equation}
where $q_m \in \mathbb{C}^{Np}$ is a column vector created by concatenating $q_k - D s_k$ vertically for $k = 1, 2, \dotsc, N$.
The function $g(h) \colon \mathbb{C} \mapsto \mathbb{C}^{Np}$ maps $h$ to the vector created by stacking the values of $C \hat{z}_i$ for $k = 1, 2, \dotsc, N$.
We remark that $g(h)$ is holomorphic and that ${\partial g (h)} / {\partial h}$ can be compiled from $C \eta_{k}$, where $\eta_{k} = {\partial \hat{z}_k} / {\partial h}$.
The latter can be computed recursively as
\begin{equation}\label{eq:css_derivetives}
\eta_{k+1} = [ (A - L C) \hat{z}_k + (B - L D) s_k + L q_k + h (A - L C) \eta_{k}] e^{-j \omega}
\end{equation}
with the initial value $\eta_{0} = 0$.
The Gauss--Newton search direction for minimizing~(\ref{eq:estObjFun_exp}) is written as follows \cite{KreutzDelgado2009}:
\begin{IEEEeqnarray}{rCl}
\delta \hat{h}^{(i)} & = & \left( \mathcal{H}^{(i)} \right)^{-1} \left( \frac{\partial g(\hat{h}^{(i)})}{\partial \hat{h}^{(i)}} \right)^H \left( q_m - g(\hat{h}^{(i)}) \right) \IEEEeqnarraynumspace \label{eq:complexGN} \\
\mathcal{H}^{(i)} & = & \left( \frac{\partial g(\hat{h}^{(i)})}{\partial \hat{h}^{(i)}} \right)^H
\left( \frac{\partial g(\hat{h}^{(i)})}{\partial \hat{h}^{(i)}} \right) \label{eq:hessianGN}
\end{IEEEeqnarray}
where $\mathcal{H}^{(i)}$ is an approximation of the Hessian matrix and the superscript $H$ indicates the conjugate transpose.

We point out that the previous optimization can be formulated outside the complex-valued framework; $J_{est,N}$ can be viewed as a function of two real-valued arguments, namely, $\Re h$ and $\Im h$.
Likewise, the estimators can be regarded as linear and real-valued, as in (\ref{eq:lti_realch_state}-\ref{eq:lti_realch_output}).
If the optimization is performed in a real-valued context, the search direction and Hessian approximation are equivalent to (\ref{eq:complexGN}-\ref{eq:hessianGN}) \cite{KreutzDelgado2009}.
In the following, we state that a function is differentiable if the derivatives of the function exist in the context of the $\mathbb{CR}$-calculus.
We briefly outline the basic relationship between the real and $\mathbb{CR}$-calculus derivatives in Appendix~\ref{app:complex_derivatives}.

Furthermore, we formulate the recursive version of the PEM (RPEM) \cite{Ljung1981}.
The update at each time step $k$ is given as
\begin{IEEEeqnarray}{rCl}
\hat{q}_{k} & = & C \hat{z}_k + D s_k \label{eq:rpem_est_output}\\
\tilde{h}_{k+1} & = & \hat{h}_k + \gamma_k \hat{S}_{k+1}^{-1} \eta_{k}^H C^H (q_k - \hat{q}_k ) \label{eq:rpem_hest} \\
\hat{S}_{k+1} & = & \hat{S}_{k} + \gamma_k [ \eta_{k}^H C^H C \eta_{k} - \hat{S}_{k} + \mu_e ] \label{eq:rpem_hess_approx} \IEEEeqnarraynumspace \\
\hat{z}_{k+1} & = & \hat{h}_{k+1} [ (A - L C) \hat{z}_k + (B - L D) s_k + L q_k  ] e^{-j \omega_k} \label{eq:rpem_state_kalman} \\
\eta_{k+1} & = & \frac{\hat{z}_{k+1}}{\hat{h}_{k+1}} + \hat{h}_{k+1} (A - L C) \eta_{k} e^{-j \omega_k}  \label{eq:rpem_state_der}
\end{IEEEeqnarray}
where $\hat{S}_{k}$ approximates $\mathcal{H}^{(k)}$ and $\mu_e \geq 0$ introduces damping to the iterative procedure.
The gain, $\gamma_k$, is a sequence of positive scalars tending to zero that weights the information contained in the current observation in relation to past observations (see also Assumption~\ref{assum:gamma_gain} for the limitation on $\gamma_{k}$).
The algorithm is initialized with the state $\hat{z}_0 = z_0$, frequency disturbance $\hat{h}_0 = h_0$ and $\hat{S}_{0} = S_0 \geq \mu_e$.

The convergence analysis of general recursive identification algorithms can be performed by associating the asymptotic behavior parameter update with an ordinary differential equation (ODE) \cite{Ljung1977,Ljung1983} and studying the stability of the ODE.
Specifically, for the RPEM, a detailed discussion on its convergence properties can be found in \cite{Ljung1981}; we transfer the results in the present case after introducing the necessary assumptions, which we discuss in the following:
\begin{assum}\label{assum:realsys_stable}
The real system is described by~(\ref{eq:linear_complex_state}-\ref{eq:linear_complex_output}) and is exponentially stable.
\end{assum}
\begin{assum}\label{assum:predictor_lti}
The model set, $\mathcal{M}$, is described by (\ref{eq:css_est_state}-\ref{eq:css_est_output}) with $h \in \mathcal{D}_\mathcal{M}$ and $\mathcal{D}_\mathcal{M} = \{ h \mid h \in \mathbb{C}, \, \cmag{h \lambda_i} \leq d_\mathcal{M}, \, i = 1, 2, \dotsc, n \}$, where $\lambda_1, \, \lambda_2, \dotsc, \lambda_n$ denotes the eigenvalues of $A - L C$. The variable $d_\mathcal{M} < 1$ is the maximum eigenvalue magnitude.
\end{assum}
\begin{assum}\label{assum:bounded_input}
The input sequence $s_k$ is bounded. The values of $s_k$ and $\omega_k$ at $k$ are independent of past values.
\end{assum}
\begin{assum}\label{assum:gamma_gain}
The sequence $\gamma_k$ satisfies $\lim_{k \to \infty} k \gamma_k = c_2$, $c_2 > 0$.
\end{assum}
For a fixed value of $h$, we denote the limit as $N \to \infty$:
\begin{equation}\label{eq:estObjFunInf}
\frac{1}{2 N} \sum_{k = 1}^{N} E [\rho_k^H \rho_k] \to J_{est,a}(h).
\end{equation}
For the sake of simplicity, we introduce the operator $\bar{E}$ acting on the function $f$ as $\bar{E}[f_k] = \lim_{N \to \infty} \frac{1}{N} \sum_{k=1}^N E[f_k]$.
Therefore, we write $ J_{est,a} = \frac{1}{2} \bar{E}[\rho_k^H \rho_k]$.
Finally, we introduce the projection method of $\tilde{h}_k$ into $\mathcal{D}_\mathcal{M}$:
\begin{equation}\label{eq:projectionM}
\hat{h}_{k+1} = \begin{cases}
\tilde{h}_{k+1}, & \tilde{h}_{k+1} \in \mathcal{D}_\mathcal{M} \\
\hat{h}_{k}, & \tilde{h}_{k+1} \notin \mathcal{D}_\mathcal{M}
\end{cases}.
\end{equation}

\begin{theorem}\label{th:openloop}
Consider the system (\ref{eq:linear_complex_state}-\ref{eq:linear_complex_output}) and Assumptions~\ref{assum:realsys_stable}, \ref{assum:predictor_lti}, \ref{assum:bounded_input} and \ref{assum:gamma_gain}.
Then, the sequence $\hat{h}_k$, which is calculated by (\ref{eq:rpem_est_output}-\ref{eq:rpem_state_der}) and the projection (\ref{eq:projectionM}), converges to a local minimum of $J_{est,a}$ or to the boundary of $\mathcal{D}_\mathcal{M}$ as $k \to \infty$.
\end{theorem}

\begin{proof}
The proof is given in Appendix~\ref{app:proofThopen}.
Here, we state general remarks on the above assumptions.
To associate the update of $\hat{h}_{k}$ with the ODE, the limit~(\ref{eq:estObjFunInf}) should be well defined.
The sequence $q_{k}$ should be bounded, which can be ensured for a stable system (Assumption~\ref{assum:realsys_stable}) and bounded input $s_{k}$ (Assumption~\ref{assum:bounded_input}).
Similar stability requirements are imposed on the estimator model to establish that both $\hat{z}_k$ and $\eta_k$ are bounded; the estimator model should be stable for constant $h$ (Assumption~\ref{assum:predictor_lti}).
\end{proof}

For the sake of completeness, we mention the associated ODEs for~(\ref{eq:rpem_hest}) and~(\ref{eq:rpem_hess_approx}):
\begin{IEEEeqnarray}{rCl}
\frac{\mathrm{d} h_D (\tau)}{\mathrm{d} \tau} & = & c_2 S_D^{-1}(\tau) f(h_D(\tau)) = -c_2 S_D^{-1}(\tau) \frac{\mathrm{d} J_{est,a}}{\mathrm{d} h_D} \label{eq:ode_f} \\
\frac{\mathrm{d} S_D (\tau)}{\mathrm{d} \tau} & = & c_2 [ F(h_D(\tau)) + \mu_e - S_D(\tau) ] \label{eq:ode_G}
\end{IEEEeqnarray}
where the subscript $D$ discriminates between the variables of the recursive algorithm and the variables of the associated ODEs.
The fictitious time $\tau$ depends on the sequence $\gamma_k$, $ f(h) = \bar{E} [ \eta_{k}^H C^H \rho_k ] $ and $ F(h) = \bar{E} [ \eta_{k}^H C^H C  \eta_{k} ]$.
Moreover, $J_{est,a}$ is the Lyapunov function used in the stability analysis of the ODEs.

Although the convergence criterion for the RPEM assumes a sequence $\gamma_k$ that tends to zero asymptotically, in practical applications, a constant value can be used if the system parameters change gradually \cite{Ljung1999}.
For systems with sudden parameter changes, a variable $\gamma_k$ scheme can be applied \cite{Fortescue1981}.
Therefore, time-varying systems can be handled by employing an appropriate choice of $\gamma_k$.
However, from a system analysis point of view, constant or adaptive $\gamma_k$ schemes are not covered by Theorem~\ref{th:openloop} and require separate analysis, see for example \cite{Anderson1986}.

\subsection{Update of the excitation frequency}\label{sec:update_freq}

The current section discusses the update of the excitation frequency, $\omega_k$, using the estimates $\hat{h}_k$.
A straightforward approach is to select $\omega_k$ as
\begin{equation}\label{eq:omega_update}
\omega_{k+1} = \arg \lambda_s + \arg \hat{h}_k = \omega_\lambda + \arg \hat{h}_k
\end{equation}
where $\omega_\lambda$ is the resonance frequency of interest corresponding to the eigenvalue $\lambda_s$ of $A$.
Determining the convergence of the RPEM under feedback (\ref{eq:omega_update}) presents two hurdles. First, although the associated ODEs presented in the asymptotic analysis of the RPEM still apply, (\ref{eq:estObjFunInf}) is not a Lyapunov function for the system because it does not account for the effect of the feedback; therefore, these ODEs cannot be used to conclude the convergence properties of the closed-loop system~\cite{Ljung1983}.
Second, the correlation between $\eta_k$ and $\epsilon_k$ further complicates the analysis.
Nonetheless, a slight modification of the algorithm can address these points.
The update of $\tilde{h}_k$ can be selected to satisfy some alternative Lyapunov function.
Constructing an estimate for $\eta_k$ that does not depend on $q_k$ removes the correlation with $\rho_k$ (referred to as the method of instrumental variables (IV)).
An uncorrelated estimate for $\eta_k$ can be realized by setting $L = 0$.
For the closed-loop system, we modify the assumptions for $\mathcal{M}$ and the input:
\begin{assum}\label{assum:predictor_iv}
The model set $\mathcal{M}$ is described by (\ref{eq:lti_ch_state}-\ref{eq:lti_ch_output}) and is restricted such that $h \in \mathcal{D}_c$, with $\mathcal{D}_c$ being a compact set. The eigenvalues of $h A$ lie strictly inside the unit circle $\forall h \in \mathcal{D}_c$.
\end{assum}
\begin{assum}\label{assum:const_input}
The input is set to a constant value, $s_k = s$.
\end{assum}
By incorporating this modification and the above assumptions, the closed-loop system becomes
\begin{IEEEeqnarray}{rCl}
q_k & = & \tilde{C} z_k + \tilde{D} s + v_k \label{eq:close_output} \\
z_{k+1} & = & (\tilde{A} z_k + \tilde{B} s + w_k) e^{-j \omega_k} \label{eq:close_system} \\
\hat{z}_{k+1} & = & \hat{h}_{k+1} (A \hat{z}_k + B s ) e^{-j \omega_k} \label{eq:close_est} \\
\eta_{k+1} & = & F \eta_k + \hat{G}_z \hat{z}_k + G_s s \label{eq:close_eta} \\
\tilde{h}_{k+1} & = & \hat{h}_k + \gamma_k \hat{S}_{k+1}^{-1} \eta_{k}^H C_\eta^H \Lambda \rho_k \label{eq:close_hest} \\
\hat{S}_{k+1} & = & \hat{S}_{k} + \gamma_k [ \eta_{k}^H C_\eta^H \Lambda C_\eta \eta_{k} - \hat{S}_{k} + \mu_e ] \label{eq:close_Sest} \IEEEeqnarraynumspace
\end{IEEEeqnarray}
where $e^{j \omega_k} = \hat{h}_k \cmag{\hat{h}_k}^{-1}$, $\Lambda$ is a constant positive definite matrix, and $C_\eta$ is a constant matrix of appropriate dimensions.
The matrices $F$, $\hat{G}_z$ and $G_s$ may depend on $\hat{h}_k$; in this case, the following restrictions apply.
\begin{assum}\label{assum:etasys_stable}
The matrices $F(h)$, $\hat{G}_z(h)$ and $G_s(h)$ are differentiable with respect to $h$, and the eigenvalues of $F(h)$ lie strictly inside the unit circle for all $h \in \mathcal{D}_c$.
\end{assum}
In the following, we do not explicitly state the dependence of $F$, $\hat{G}_z$ and $G_s$ on $h$; this dependence should be assumed unless stated otherwise.

\begin{theorem}\label{th:closed_general}
Consider the system (\ref{eq:close_output}-\ref{eq:close_Sest}) and Assumptions~\ref{assum:realsys_stable}, \ref{assum:gamma_gain}, \ref{assum:predictor_iv}, \ref{assum:const_input} and \ref{assum:etasys_stable}, along with a projection that always maintains $\hat{h}_k \in \mathcal{D}_c$.
Assume that there exists a real positive function $V_D(h_D,S_D)$ such that
\begin{equation}
\frac{\mathrm{d} V_D (h_D(\tau),S_D(\tau))}{\mathrm{d} \tau} < 0, \quad h_D \in \mathcal{D}_c
\end{equation}
along the trajectories
\begin{IEEEeqnarray}{rCl}
\frac{\mathrm{d} h_D (\tau)}{\mathrm{d} \tau} & = & c_2 S_D^{-1}(\tau) \eta_a^H C_\eta^H \Lambda \rho_a \label{eq:ode_f_closed} \\
\frac{\mathrm{d} S_D (\tau)}{\mathrm{d} \tau} & = & c_2 [ \eta_a^H C_\eta^H \Lambda C_\eta \eta_a + \mu_e - S_D(\tau) ] \label{eq:ode_G_closed}
\end{IEEEeqnarray}
with $\eta_a = (I-F)^{-1}(\hat{G}_z \hat{z}_a + G_s s)$, $\rho_a = \tilde{C} z_a - C \hat{z}_a + (\tilde{D} - D) s$, $\hat{z}_a = h_D \hat{H}^{-1} B s $ and $z_a = \tilde{H}^{-1} \tilde{B} s $.
The matrices $\tilde{H}$ and $\hat{H}$ are defined as $\tilde{H} = I e^{j \omega(h_D)} - \tilde{A}$ and $\hat{H} = I e^{j \omega(h_D)} - h_D A$, respectively.
Let
\begin{equation}
\mathcal{D}_V = \{ h, S \mid \frac{\mathrm{d} V_D (h,S)}{\mathrm{d} \tau} = 0\}.   
\end{equation}
Then, as $k \to \infty$, $\{ \hat{h}_k, \hat{S}_k \}$ tends to $\mathcal{D}_V$, or $\hat{h}_k$ tends to the boundary of $\mathcal{D}_c$.
\end{theorem}

\begin{proof}
The asymptotic analysis of (\ref{eq:close_output}-\ref{eq:close_Sest}) is based on \cite[theorem~4.2]{Ljung1983}, which applies the technique of the associated ODEs to a general family of recursive algorithms.
We verify the necessary conditions for the application of the theorem (labeled Conditions M1, M2, Cr1, Cr2, R1, G1 and A1 in \cite{Ljung1983}) and confirm the requirement that the system is described by a linear structure.

The closed-loop system can be written as
\begin{IEEEeqnarray}{rCl}
z_{c,k+1} & = & \mathcal{A}(h) z_{c,k} + \mathcal{B}_s (h) s + \mathcal{B}_w w_k \label{eq:closed_ss_state} \\
\begin{bmatrix}
\rho_k \\ \eta_{k}
\end{bmatrix} & = &
\begin{bmatrix}
\mathcal{C}_1 \\ \mathcal{C}_2
\end{bmatrix} z_{c,k} +
\begin{bmatrix}
\tilde{D} - D \\ 0
\end{bmatrix} s +
\begin{bmatrix}
I \\ 0
\end{bmatrix} v_k \label{eq:closed_ss_output}
\end{IEEEeqnarray}
where $z_{c,k} = [z_k^T, \hat{z}_k^T, \eta_{k}^T]^T$.
Taking into account that $e^{j \omega(h)} = h \cmag{h}^{-1}$, the system matrices are given as
\begin{IEEEeqnarray}{rCl}
\mathcal{A}(h) & = &
\begin{bmatrix}
\overline{h} \cmag{h}^{-1} \tilde{A} & 0 & 0  \\
0 &  \cmag{h} A & 0 \\
0 & \hat{G}_z & F
\end{bmatrix} \IEEEnonumber \\
\mathcal{B}_s(h) & = &
\begin{bmatrix}
\overline{h} \cmag{h}^{-1} \tilde{B}  \\
\cmag{h} B \\
G_s
\end{bmatrix} ,\quad
\mathcal{B}_w =
\begin{bmatrix}
I \\
0 \\
0
\end{bmatrix} \IEEEnonumber \\
\begin{bmatrix}
\mathcal{C}_1 \\ \mathcal{C}_2
\end{bmatrix} & = &
\begin{bmatrix}
\tilde{C} & -C & 0 \\
0 & \phantom{-}0 & I
\end{bmatrix} \IEEEnonumber
\end{IEEEeqnarray}
which constitutes a linear model structure.
Assumptions~\ref{assum:realsys_stable}, \ref{assum:predictor_iv} and \ref{assum:etasys_stable} ensure that (\ref{eq:closed_ss_state}-\ref{eq:closed_ss_output}) is stable for all $h \in \mathcal{D}_c$ and is differentiable with respect to $h$, satisfying Conditions~M1 and M2.

Condition~Cr1 sets the smoothness requirements for the function that determines the update for $\hat{h}_k$, $\tilde{f}(k, h, \rho, \eta) = \eta^H C_\eta^H \Lambda \rho$.
Since the conditions in \cite{Ljung1983} assume real-valued functions, we assume that the norms are taken as if the functions are real-valued; we have adapted the relations to take into account the differences that arise from the complex nature of the problem, as described in Appendix~\ref{app:complex_derivatives}.
The function $\tilde{f}$ is differentiable with respect to $h$, $\rho$ and $\eta$, and
\begin{equation}
\lVert \tilde{f} \rVert + \left \lVert \frac{\partial \tilde{f}}{ \partial h} \right \rVert = \lVert \eta^H C_\eta^H \Lambda \rho \rVert \leq \lVert C_\eta \eta  \rVert \lVert \Lambda \rho \rVert \leq \frac{1}{2}(\lVert \Lambda \rho \rVert^2 + \lVert C_\eta \eta  \rVert^2)  \leq c_3 ( 1 + \lVert \rho \rVert^2 + \lVert \eta \rVert^2).
\end{equation}
for some $c_3 < \infty$. In the previous derivation, we used the Cauchy--Schwarz inequality.
Likewise,
\begin{equation}
\left \lVert \frac{\partial \tilde{f}}{ \partial \rho} \right \rVert + \left \lVert \frac{\partial \tilde{f}}{ \partial \eta } \right \rVert = \lVert \eta^H C_\eta \Lambda \rVert + \lVert C_\eta \Lambda \rho \rVert \leq c_3 ( 1 + \lVert \rho \rVert + \lVert \eta \rVert).
\end{equation}
Similar smoothness conditions must be verified for the update function of $\hat{S}_k$, $\tilde{F}(k,S,h,\rho,\eta) = \eta^H C_\eta^H \Lambda C_\eta \eta - S$.
The function should be differentiable with respect to $S$, $h$, $\rho$ and $\eta$, which is true, and
\begin{IEEEeqnarray}{rCl}
\lVert \tilde{F} \rVert + \left \lVert \frac{\partial \tilde{F}}{ \partial S} \right \rVert + \left \lVert \frac{\partial \tilde{F}}{ \partial h} \right \rVert & = & \lVert \eta^H C_\eta^H \Lambda C_\eta \eta - S \rVert + 1 \leq c_3 (1 + \lVert \rho \rVert^2 + \lVert \eta \rVert^2) \\
\left \lVert \frac{\partial \tilde{F}}{ \partial \rho} \right \rVert + \left \lVert \frac{\partial \tilde{F}}{ \partial \eta} \right \rVert & \leq & 2 \lVert \eta^H C_\eta^H \Lambda C_\eta \rVert \leq c_3 (1 + \lVert \rho \rVert + \lVert \eta \rVert).
\end{IEEEeqnarray}
Therefore, $\tilde{F}$ complies with Condition~Cr2.

We presume that $\hat{S}_k \geq \mu_e > 0$; then,
\begin{equation}
\hat{S}_{k+1} = (1-\gamma_k) \hat{S}_k + \gamma_k (\eta_{k}^H C_\eta^H \Lambda C_\eta \eta_{k} + \mu_e) \geq (1-\gamma_k) \hat{S}_k + \gamma_k \mu_e \geq \mu_e
\end{equation}
since $\eta_{k}^H C_\eta^H \Lambda C_\eta \eta_k \geq 0$.
Therefore, if $\hat{S}_0 \geq \mu_e$, then by mathematical induction, $\hat{S}_k \geq \mu_e, \, \forall k > 0$, which is in accordance with Condition~R1.
Condition~G1 coincides with Assumption~\ref{assum:gamma_gain}.

Condition~A1 sets two requirements. First, the input sequences $s$, $w_k$ and $v_k$ are such that $f(h) = \bar{E}[\eta_k^H C_\eta^H \Lambda \rho_k]$ and $F(h) = \bar{E}[\eta_k^H C_\eta^H \Lambda C_\eta \eta_k]$ exist $\forall h \in \mathcal{D}_c$.
The limits are well defined since $s$ is bounded, $w_k$ and $v_k$ are stationary, and $\rho_k$ is asymptotically stationary.
Moreover, $w_k$ and $v_k$ have bounded moments and
\begin{equation}
\lim_{N \to \infty} \sup \frac{1}{N} \sum_{k=1}^{N}(1 + \lVert s \rVert + \lVert w_k \rVert + \lVert v_k \rVert)^3 < \infty
\end{equation}
satisfying the second requirement of Condition~A1.
The aforementioned limits can be written in closed form because for a fixed $h$, (\ref{eq:closed_ss_state}) approaches a steady state.
We denote $z_a = \bar{E}[z_k]$, $\hat{z}_a = \bar{E}[\hat{z}_k]$, $\rho_a = \bar{E}[\rho_k]$ and $\eta_a = \bar{E}[\eta_k]$, which are computed as
\begin{IEEEeqnarray}{rCl}
z_a & = & (I e^{j \omega(h)} - \tilde{A})^{-1} \tilde{B} s = \tilde{H}^{-1} \tilde{B} s \\
\hat{z}_a & = & (I e^{j \omega(h)} - h A)^{-1} h B s = \hat{H}^{-1} h B s \\
\rho_a & = & \tilde{C} z_a - C \hat{z}_a + (\tilde{D} - D) s \\
\eta_a & = & (I - F)^{-1} (\hat{G}_z \hat{z}_a + G_s s).
\end{IEEEeqnarray}
Since $\eta_k$ is not correlated with $\rho_k$, $f(h) = \eta_a^H C_\eta^H \Lambda \rho_a$ and $F(h) = \eta_a^H C_\eta^H \Lambda C_\eta \eta_a$.

All of the conditions for the application of \cite[theorem~4.1]{Ljung1983} are satisfied, and the associated ODEs are given in (\ref{eq:ode_f_closed}-\ref{eq:ode_G_closed}).
Therefore, if the function $V_D$ exists and is strictly decreasing along the trajectories of (\ref{eq:ode_f_closed}-\ref{eq:ode_G_closed}), then $\{\hat{h}_k, \hat{S}_k \}$ tends to $\mathcal{D}_V$ or $\hat{h}_k$, the boundary of $\mathcal{D}_c$, as $k \to \infty$.
\end{proof}

Selecting a general Lyapunov function is a challenging task.
The process can be simplified if we assume that the real system is contained in $\mathcal{M}$, written formally as
\begin{IEEEeqnarray}{rCl}
z_{k+1} & = & h_s (A z_k + B s + w_k) e^{j \omega(h)} \label{eq:sysinm_system} \\
q_{k} & = & C z_k + D s + v_k \label{eq:sysinm_output}
\end{IEEEeqnarray}
with $h_s \in \mathcal{D}_c$.
In this case, the prediction error for the state, $\epsilon_k = z_k - \hat{z}_k$, evolves as
\begin{equation}
\epsilon_{k+1} = [h_s (A \epsilon_k + w_k) + (h_s - \hat{h}_k)(A \hat{z}_k + B s) ] e^{-j \omega_k}. \>
\end{equation}
\begin{theorem}\label{th:closed_inm}
Consider the system (\ref{eq:sysinm_system}-\ref{eq:sysinm_output}) and Assumptions~\ref{assum:gamma_gain}, \ref{assum:predictor_iv}, \ref{assum:const_input} and \ref{assum:etasys_stable}. Assume a projection that always maintains $\hat{h}_k \in \mathcal{D}_c$.
Then, the sequence of $\hat{h}_k$ produced by the algorithm (\ref{eq:close_est}-\ref{eq:close_Sest}) converges asymptotically to $h_s$ or to the boundary of $\mathcal{D}_c$ if $M + M^H > 0, \, \forall h \in \mathcal{D}_c$, with
\begin{equation}
M(h) = s^H [C_\eta (I - F)^{-1} (h \hat{G}_z \hat{H}^{-1} B+ G_s)]^H \Lambda C H^{-1} \hat{H}^{-1} B s e^{j \omega(h)} \label{eq:closed_inm_matrix}
\end{equation}
and $H = I e^{j \omega(h)} - h_s A$.
\end{theorem}
\begin{proof}
We apply Theorem~\ref{th:closed_general} in the case that $\tilde{A} = h_s A$, $\tilde{B} = h_s B$, $\tilde{C} = C$ and $\tilde{D} = D$.
The candidate Lyapunov function is $V_D = \overline{(h_s - h)} (h_s - h)$.
The asymptotic values for $z_a$ and $\epsilon_a = \bar{E}[\epsilon_k]$ become
\begin{IEEEeqnarray}{rCl}
z_a & = & (I e^{j \omega(h)} - h_s A)^{-1} h_s B s = H^{-1} h_s B s \\
\epsilon_a & = & (I e^{j \omega(h)} - h_s A)^{-1} \frac{\hat{z}_a}{h} (h_s - h) e^{j \omega(h)} = H^{-1} \frac{\hat{z}_a}{h} (h_s - h) e^{j \omega(h)}
\end{IEEEeqnarray}
with an estimation error of $\rho_a = C \epsilon_a$.
The rate of change of $V_D$ along the trajectory of (\ref{eq:ode_f_closed}-\ref{eq:ode_G_closed}) is given as
\begin{IEEEeqnarray}{rCl}
\frac{\mathrm{d} V_D}{\mathrm{d} \tau} & = & \frac{\partial V_D}{\partial h_D} \frac{\mathrm{d} h_D}{\mathrm{d} \tau} + \frac{\partial V_D}{\partial \overline{h_D}} \frac{\mathrm{d} \overline{h_D}} {\mathrm{d} \tau} = -(M + M^H) \frac{c_2 V_D}{S_D} \IEEEeqnarraynumspace \\
M & = & [C_\eta (I - F)^{-1} (\hat{G}_z \hat{z}_a + G_s s)]^H \Lambda C H^{-1} \frac{\hat{z}_a}{h} e^{j \omega(h)}. \IEEEnonumber
\end{IEEEeqnarray}
The Lyapunov function is a decreasing function if $ M + M^H > 0$ since $V_D$, $c_2$ and $S_D$ are positive real variables.
Starting from some initial $\hat{h}_0 \in \mathcal{D}_c$, as $k \to \infty$, the algorithm (\ref{eq:close_est}-\ref{eq:close_Sest}) drives $\hat{h}_k \to h_s$ or the boundary of $\mathcal{D}_c$.
\end{proof}

Theorem~\ref{th:closed_inm} allows one either to select $F$, $\hat{G}_z$, $G_s$, $C_\eta$ and $\Lambda$ and ensure the convergence of $h$ to $h_s$ in a given $\mathcal{D}_c$ or to estimate $\mathcal{D}_c$ for a given update of $\eta_k$.
For models with perfect state information, the analysis of the resonance tracking algorithm can be significantly simplified.
We introduce the following lemma to facilitate the discussion.
\begin{lem}\label{lem:mat_positive_def}
Given the invertible matrices $X, Y \in \mathbb{C}^{n \times n}$ and $Z = Y - X$,
\begin{IEEEeqnarray}{l}
Y Y^H - Z Z^H > 0 \Leftrightarrow \IEEEnonumber \\
(Y^{-1})^H X^{-1} + (X^{-1})^H Y^{-1} - (Y^{-1})^H Y^{-1} > 0. \IEEEeqnarraynumspace
\end{IEEEeqnarray}
\end{lem}
\begin{proof}
We substitute the expansion for $Z$ into $Y Y^H - Z Z^H \ge 0$:
\begin{IEEEeqnarray*}{rCl}
Y Y^H - (Y - X) (Y - X)^H & > & 0 \Leftrightarrow \\
Y X^H + X Y^H - X X^H  & > & 0 \Leftrightarrow \\
X^{-1} Y + Y^H (X^{-1})^H - I & > & 0 \Leftrightarrow \\
(Y^{-1})^H X^{-1} + (X^{-1})^H Y^{-1} - (Y^{-1})^H Y^{-1} & > & 0.
\end{IEEEeqnarray*}
\end{proof}
\begin{corol}
Given the system (\ref{eq:sysinm_system}-\ref{eq:sysinm_output}) with $C \in \mathbb{R}^{n \times n}$ being invertible, selecting $F = 0$, $\hat{G}_z = A e^{-j \omega_k}$, $G_s = B e^{-j \omega_k}$, $C_\eta = (C^{-1})^H$ and $\Lambda = I$ ensures that for $\mathcal{D}_c = \{ h \mid h \in \mathbb{C}, \, \lVert h A \rVert \leq d_M < 1 \}$, the algorithm (\ref{eq:close_est}-\ref{eq:close_Sest}) will converge to $h_s$ or to the boundary of $\mathcal{D}_c$ with an appropriate projection for $\hat{h}_k \in \mathcal{D}_c$.
\end{corol}
\begin{proof}
The previous statement is a special case of Theorem~\ref{th:closed_inm}.
The asymptotic value for $\hat{z}_a$ satisfies $h^{-1} \hat{z}_a = (A \hat{z}_a +  B s) e^{-j \omega(h)}$, which results in $ M = \hat{z}_a^H H^{-1} \hat{z}_a \cmag{h}^{-2} e^{j \omega(h)}$.
We apply Lemma~\ref{lem:mat_positive_def} with $X = H$ and $Y = I e^{j \omega(h)}$.
We note that $Z = I e^{j \omega(h)} - H = h_s A$.
From the definition of $\mathcal{D}_c$, we have $\lVert h_s A \rVert < 1$, which implies that $I - \cmag{h_s}^2 A^H A > 0$ and $e^{j \omega(h)} H^{-1} + e^{-j \omega(h)} (H^{-1})^H - I > 0$.
Thus, $M + M^H > 0$ for $\rVert \hat{z}_a \rVert > 0$, concluding the proof.
\end{proof}

It is more relevant for the imperfect state information setting to assign the update for $\eta_k$ and to identify the set $\mathcal{D}_c$ that establishes convergence to $h_s$.
The estimation of $\mathcal{D}_c$ can be formulated as a constraint optimization problem.
For simplicity, we parameterize $\mathcal{D}_c$ as an annular sector on the complex plane $\mathcal{D}_c = \{ h \mid h \in \mathbb{C}, \, d_m \leq \cmag{h} \leq d_M, \vert \arg h \vert \leq d_\phi \}$.
The maximum radius $d_M$ can be set to a value close to $\lVert A \rVert^{-1}$, satisfying Assumption~\ref{assum:predictor_iv}.
We seek to maximize the area of the sector:
\begin{IEEEeqnarray}{rCl}
\max_{d_m, d_\phi} & \quad & d_\phi (d_M^2 - d_m^2) \label{eq:optim_LyapReg} \\
\text{s.t.} & & \min_{h, h_s \in \mathcal{D}_c} (M + M^H) > 0 \IEEEnonumber \\
& & \vert d_\phi \vert \leq \pi \IEEEnonumber  \\
& & 0 \leq d_m < d_M. \IEEEnonumber
\end{IEEEeqnarray}
Taking into account (\ref{eq:omega_update}), the dimensionality of the minimization of the first constraint can be reduced for some updates.
Specifically, $H e^{-j \omega(h)} = I - \cmag{h} e^{-j \omega_\lambda} A$ and $\hat{H} e^{-j \omega(h)} = I - \cmag{h_s} e^{j (\delta \phi - \omega_\lambda)} A$, where $\delta \phi = \arg h_s - \arg h$.
The constraint can be confirmed by minimizing $M + M^H$ with respect to $d_m \leq \cmag{h} \leq d_M$, $d_m \leq \cmag{h_s} \leq d_M$ and $ \vert \delta \phi \vert \leq 2 d_\phi$.
We note that $M+M^H \in \mathbb{R}$ and the problem can be solved using standard optimization algorithms.

Finally, we study a noise-free system.
In the absence of disturbances and measurement noise, there is no correlation between $\eta_k$ and $\epsilon_k$, and a linear term of $q_k$ can be included in the update for $\eta_k$.
The noise-free system equations and the modified $\eta_k$ update become
\begin{IEEEeqnarray}{rCl}
q_{k} & = & C z_k + D s \label{eq:noisefree_output} \\
z_{k+1} & = & h_s (A z_k + B s) e^{j \omega(h)} \label{eq:noisefree_system} \\
\eta_{k+1} & = & F \eta_k + G_q q_k + \hat{G}_z \hat{z}_k + G_s s\label{eq:noisefree_eta}.
\end{IEEEeqnarray}
\begin{theorem}\label{th:nonoise}
Consider the system (\ref{eq:noisefree_output}-\ref{eq:noisefree_system}) and Assumptions~\ref{assum:gamma_gain}, \ref{assum:predictor_iv}, \ref{assum:const_input} and \ref{assum:etasys_stable}. Assume a projection that always maintains $\hat{h}_k \in \mathcal{D}_c$.
Then, the sequence of $\hat{h}_k$ produced by the algorithm (\ref{eq:close_est}), (\ref{eq:noisefree_eta}) and (\ref{eq:close_hest}-\ref{eq:close_Sest}) converges asymptotically to $h_s$ or to the boundary of $\mathcal{D}_c$ if $M + M^H > 0, \, \forall h \in \mathcal{D}_c$, with
\begin{equation}
M(h)= s^H \big \{ C_\eta (I - F)^{-1} [G_q (h_s C H^{-1} B + D) + G_s + h \hat{G}_z \hat{H}^{-1} B ]  \big \}^H \Lambda C H^{-1} \hat{H}^{-1} B s e^{j \omega(h)}. \label{eq:nonoise_matrix}
\end{equation}
\end{theorem}
\begin{proof}
The proof is analogous to Theorem~\ref{th:closed_inm}.
\end{proof}

For a system with perfect state information, similar results can be stated for the noise-free case.
\begin{corol}
Given a system (\ref{eq:noisefree_output}-\ref{eq:noisefree_system}), with $C \in \mathbb{R}^{n \times n}$ being invertible,
selecting $F = \hat{G}_z = 0$, $G_q = A C^{-1} e^{j \omega(h)}$, $G_s = (B - A C^{-1} D) e^{j \omega(h)}$ and $C_\eta = (C^{-1})^H$ ensures that for $\mathcal{D}_c = \{ h \mid h \in \mathbb{C}, \, \lVert h A \rVert \leq d_M < 1 \}$, the algorithm (\ref{eq:close_est}), (\ref{eq:noisefree_eta}) and (\ref{eq:close_hest}-\ref{eq:close_Sest}) will converge to $h_s$ or to the boundary of $\mathcal{D}_c$ with an appropriate projection for $\hat{h}_k$.
\end{corol}
\begin{proof}
The previous statement is a special case of Theorem~\ref{th:nonoise}.
The asymptotic value for $z_a$ satisfies $h_s^{-1} z_a = (A z_a +  B s) e^{-j \omega(h)}$, and the state error can be alternatively written as
\begin{equation*}
\epsilon_a = \hat{H}^{-1} \frac{z_a}{h_s} (h_s - h) e^{j \omega(h)}
\end{equation*}
which results in $ M = z_a^H \hat{H}^{-1} z_a \cmag{h_s}^{-2} e^{j \omega(h)}$.
We apply Lemma~\ref{lem:mat_positive_def} with $X = \hat{H}$ and $Y = I e^{j \omega_k}$.
We note that $Z = I e^{j \omega(h)} - \hat{H} = h A$.
From the definition of $\mathcal{D}_c$, we have $\lVert h A \rVert < 1$, implying that $I - \cmag{h}^2 A^H A > 0$ and $e^{j \omega(h)} \hat{H}^{-1} + e^{-j \omega(h)} (\hat{H}^{-1})^H - I > 0$.
Thus, $M + M^H > 0$ for $\rVert z_a \rVert > 0$, concluding the proof.
\end{proof}

\section{Implementation and numerical simulations}\label{sec:simulations}
The present section discusses ways to estimate the complex envelope of signals in real time and evaluates the performance of the proposed algorithms through simulation.
Regarding the real-time requirements of the general tracking algorithm, each time step involves the evaluation of (\ref{eq:close_est}-\ref{eq:close_Sest}), and
the memory requirements scale with $n$.
The processing power can range significantly depending on the value of $\omega_{\lambda}$, which imposes the range of $T_{s}$ and the time interval on computations between updates $\omega_{k}$.
The number of operations per frequency update is determined by $n$.

\subsection{Estimating the analytic representation}
The main issue related to the implementation of the proposed frequency tracking scheme is the extraction of the complex envelope $q_k$ from the measured output.
The sliding discrete Fourier transform (sDFT) and sliding Goertzel algorithms have been proposed to convert a signal into its analytic representation \cite{Jacobsen2003}.
The sDFT is equivalent to the discrete Fourier transform (DFT) applied to a window of length $N_f \in \mathbb{N}_{>0}$; the output rate is equal to the input signal rate, but the analytic representation is computed only at a specified center frequency.

For the present application, we need to estimate $q_k$ at noninteger multiples of $T_s^{-1}$.
Moreover, $\omega_k$ may change over the computation window, which makes the nonuniform discrete time Fourier transform a more appropriate choice than the DFT \cite{Duijndam1999}.
Hence, we extend the sDFT algorithm to the sliding nonuniform discrete time Fourier transform (sNDTFT) case.
Following the derivation of the sDFT in \cite{Jacobsen2003}, the sNDTFT filter is formulated as
\begin{IEEEeqnarray}{rCl}
\tilde{Y}_k & = &  \tilde{Y}_{k-1} e^{j \omega_k} - y_{k-N_f} e^{j \delta \omega_k} + y_k \label{eq:sdft_state}\\
\delta \omega_k & = & \delta \omega_{k-1} + \omega_k - \omega_{k-N_f} \label{eq:sdft_binRemAngle}
\end{IEEEeqnarray}
where $ \tilde{Y}_k$ is an internal state and $\omega_k = 0$ and $y_k = 0$ for $k < 0$.
The complex envelope $q_k$ is reconstructed as
\begin{equation}\label{eq:sdft_output}
q_k = 2 \frac{ \tilde{Y}_k}{N_f} e^{-j \theta_k}
\end{equation}
where the last term has a dual role: to apply the phase correction introduced in \cite{Sysel2012}, which accounts for the calculation at a noninteger multiple of $T_s^{-1}$, and to match the phase $q_k$ with~(\ref{eq:theta}).
More details are given in Appendix~\ref{app:sDFT_der}.

To mitigate the effect of spectral leakage, a Hann window is applied in the frequency domain; we compute~(\ref{eq:sdft_state}) for two adjacent frequencies $\omega_k \pm 2\pi/N_f$, and the results are averaged and subtracted from $\tilde{Y}_k$ before the calculation~(\ref{eq:sdft_output}).
This calculation includes the correction factor of 2 needed to recover the correct signal amplitude.
Since the calculations are performed with an offset of $2\pi/N_f$, $\delta \omega_k$ remains the same for the adjacent frequencies; only the first $N_f$ samples will exhibit a mismatch.
We point out that the sNDTFT will produce an approximation of the complex envelope, which depends on $N_f$.
Additionally, the responsiveness of the algorithm to changes in $\omega_k$ is also affected by $N_f$.
The real-time implementation of the sNDTFT requires storing $N_{f}$ values per system output and an additional $N_{f}$ past values of $\omega_{k}$.

\subsection{Numerical simulations}
For the simulations, we selected systems with imperfect state information.
We tested three updates for $\eta_k$ influenced by the analysis of Section~\ref{sec:update_freq}:
\begin{itemize}
\item $\eta_{k+1} = (hA \eta_k + A \hat{z}_k + B s) e^{-j \omega_k}$, $C_\eta = C$ and $\Lambda = I$. Substitution in (\ref{eq:closed_inm_matrix}) results in $M = \hat{z}_a^H (\hat{H}^{-1})^H C^H C H^{-1} \hat{z}_a \cmag{h}^{-2}$. Due to the similarity of the present algorithm to the RPEM algorithm introduced for the open-loop system, we refer to this update as the cRPEM.
\item $\eta_k = (A \hat{z}_k + B s) e^{-j \omega_k}$, $C_\eta = c_5 C$ and $\Lambda=I$. Eq.~(\ref{eq:closed_inm_matrix}) becomes $M = c_5 \hat{z}_a^H C^H C H^{-1} \hat{z}_a \cmag{h}^{-2} e^{j \omega_k}$. We refer to this update as recursive IV (RIV).
\item  $\eta_k = q_k - D s$, $C_\eta = c_5 I$ and $\Lambda=I$. Eq.~(\ref{eq:nonoise_matrix}) becomes $M = c_5 z_a^H C^H C \hat{H}^{-1} z_a h_s^{-1} e^{j \omega_k}$. We refer to this update as the output association (OA).
\end{itemize}

We assessed the robustness of each update by calculating $\mathcal{D}_c$ through the optimization (\ref{eq:optim_LyapReg}).
We note that $\mathcal{D}_c$ is calculated under the assumption that the real system can be described by (\ref{eq:sysinm_system}-\ref{eq:sysinm_output}) using the matrices $M$ given above, which is not the case.
We used the area of $\mathcal{D}_c$ as a metric to assess the robustness of each algorithm.
Moreover, for the simulations, we used the estimated $\mathcal{D}_c$ for the projection scheme despite the fact that the simulations were performed with the system described by (\ref{eq:linear_uncertain_state}-\ref{eq:linear_uncertain_output}).
Furthermore, we implemented a projection scheme that is more appropriate for an annular sector that selects the phase and the magnitude separately:
\begin{IEEEeqnarray}{rCl}
\arg \hat{h}_{k+1} & = & \begin{cases}
-d_\phi, & \arg \tilde{h}_{k+1} < -d_\phi \\
d_\phi, & \arg \tilde{h}_{k+1} > d_\phi \\
\arg \tilde{h}_{k+1}, & \text{otherwise}
\end{cases} \label{eq:projection_phase} \\
\cmag{\hat{h}_{k+1}} & = & \begin{cases}
d_m, & \cmag{\tilde{h}_{k+1}} < d_m \\
d_M, & \cmag{\tilde{h}_{k+1}} > d_M \\
\cmag{\tilde{h}_{k+1}}, & \text{otherwise}
\end{cases} \label{eq:projection_magnitude}.
\end{IEEEeqnarray}
The numerical simulations and optimizations were performed using Scilab/Xcos 6.0.2 software.

\begin{figure}
\centering
\includegraphics{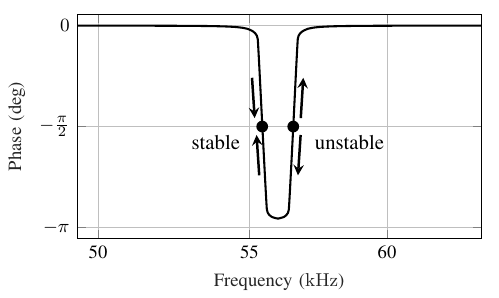}
\caption{Phase difference of the piezoelectric transducer~(\ref{eq:piezo_laplace}). The points indicate the resonance and the antiresonance at $-\pi/2$, which coincide with the equilibria for control techniques based solely on phase information (i.e., the PLL). The arrows show the attraction (stable) and repulsion (unstable) in the region around the equilibria.}
\label{fig:piezoPhase}
\end{figure}

\subsubsection{SISO system with a nonmonotonic phase}\label{sec:subsec_piezo}
Initially, we simulated the resonance tracking algorithm for the model of a piezoelectric actuator, a system with both resonance and antiresonance.
The actuator is modeled as a resistor-inductor-capacitor oscillator ($R_m$, $L_m$, and $C_m$) in parallel with a capacitor $C_0$, as described in \cite{Liu2015}.
The transfer function $G_p(\mathrm{s})$ from voltage to charge in the Laplace domain is given as
\begin{equation}\label{eq:piezo_laplace}
G_p(\mathrm{s}) = C_0 \frac{\mathrm{s}^2 + \mathrm{s} \frac{R_m}{L_m} + \frac{C_m+C_0}{L_m C_m C_0}}{\mathrm{s}^2 + \mathrm{s} \frac{R_m}{L_m} + \frac{1}{L_m C_m}}
\end{equation}
where $\mathrm{s}$ is the Laplace variable.
The nominal values for the parameters ($C_0 = \SI{2}{\nano\farad}$, $R_m = \SI{50}{\ohm}$, $L_m = \SI{0.103}{\henry}$, and $C_m = \SI{80}{\pico\farad}$) are taken from \cite{Liu2015}.

The model exhibits a nonmonotonic input-output phase, as shown in Fig.~\ref{fig:piezoPhase}, which is challenging for resonance tracking techniques based solely on phase information.
Since the phase is not unique, multiple equilibria arise, which alternate between stable and unstable modes \cite{Hayashi1992,Zhang2015,Brack2016a,Liu2015}.
The points in Fig.~\ref{fig:piezoPhase} mark the two equilibria at a phase difference of $-\pi/2$.
The stability of the equilibria for a negative gain controller is marked by arrows.
To demonstrate that the proposed resonance tracking algorithm remains unaffected by the nonmonotonic nature of the phase, we assume that the model parameters vary by 10\% around their nominal values.

\begin{figure}
\centering
\includegraphics{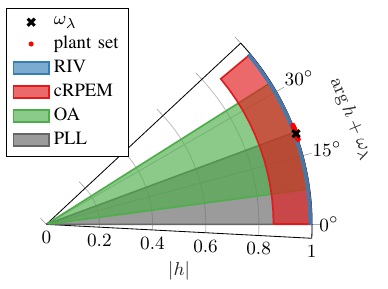}
\caption{The region $\mathcal{D}_c$ for the piezoelectric transducer~(\ref{eq:piezo_laplace}) and the different update algorithms. The region of stability for the PLL along with the $\tilde{\omega}_{\lambda}$ for a random plant set are shown for comparison.}
\label{fig:piezoPolar}
\end{figure}

First, we determined the region $\mathcal{D}_c$, where $M$ satisfies Theorems~\ref{th:closed_inm} and \ref{th:nonoise}, for the different $\eta$ updates.
The resulting annuli are plotted in Fig~\ref{fig:piezoPolar}.
The OA update covers more area on the complex plane, followed by the cRPEM and RIV.
The frequency response of (\ref{eq:piezo_laplace}) away from resonance is almost constant, limiting the range of $\mathcal{D}_c$ around $\omega_{\lambda}$.
The region of stability for the PLL algorithm, which includes all resonance shifts that do not cross the antiresonance, can also be seen in Fig~\ref{fig:piezoPolar}.
The resonances of a random set of plants with 10\% parametric uncertainity are included for comparison.

Next, we performed Monte Carlo sampling of the model parameters to create a set of one hundred random plants.
The systems were discretized with $T_s = \SI{1}{\micro\second}$ and converted into the minimal and balanced state-space realization, and we set $s = \SI{1}{\volt}$.
The noise was selected to have a power of approximately 10\% of the signal power at resonance; specifically, $Q_s = \SI{0.01}{\volt\tothe{2}}$ and $R = \SI{64}{\nano\coulomb\tothe{2}}$.
We point out that the signal-to-noise ratio will be significantly higher than 10\% ``far'' from the resonance frequency.

\begin{figure}
\centering
\includegraphics{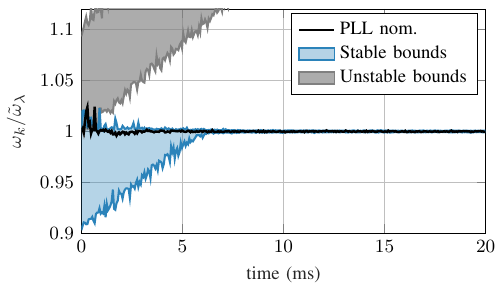}
\caption{Simulations results of the PLL algorithm with a set of 100 plants with 10\% parametric uncertainty.
The regions bounding the stable and unstable trajectories are shown along with the trajectory for the nominal plant.}
\label{fig:piezoRes_PLL}
\end{figure}

As a point of reference, we first simulated the model set with a simple PLL scheme, which we implemented as a discrete-time proportional--integral--derivative (PID) controller with a derivative filter:
\begin{equation}
C_{\rm PID}({\rm z}) = K_{p} + K_{i} T_{s} \frac{{\rm z} + 1}{{\rm z} - 1} + \frac{K_{d}}{T_{f} + \frac{T_{s} {\rm z}}{{\rm z} + 1}}
\end{equation}
where $K_{p}$, $K_{i}$ and $K_{d}$ are the proportional, integral and derivative gains, respectively, $T_{f}$ is the derivative filter time constant, and $\rm z$ is the discrete shift operator.
The phase error from the setpoint $-\pi/2$, $\delta \phi_{k} = \pi/2 + (\varphi_{k} - \psi_{k})$, is fed to $C_{\rm PID}$, and the output is added to the PLL center frequency, $\omega_{\rm PLL}$, to determine $\omega_{k}$.
We used the sNDTFT with $N_{f} = 32$ to determine $\delta \phi_{k}$, and we selected $\omega_{\rm PLL} = \omega_{\lambda}$.
The PID parameters were selected as a compromise between the responsiveness and the rejection of noise and set to $K_{p} = 0.001$, $K_{i} = 2$, $K_{d} = 0.001$ and $T_{f} = 2$.
The simulation results are shown in Fig.~\ref{fig:piezoRes_PLL}, showing that for 47\% of the uncertain plants, the PLL failed to lock on $\tilde{\omega}_{\lambda}$.
This result is in accordance with Fig.~\ref{fig:piezoPolar}, where almost half of the uncertain plants lie outside the stable PLL region.
We note that by setting $\omega_{\rm PLL}$ to a significantly lower value than $\omega_{\lambda}$, the region of instability could have been avoided during the initial transient at the cost of a slower response.
This approach does not avoid disturbances that can perturb the controller from the resonance lock into the region of instability.

For the cRPEM, we set $\gamma = 0.002$, $\hat{S}_{0} = \SI{0.01}{\micro\coulomb\tothe{2}}$ and $\mu_e = \SI{8e-4}{\micro\coulomb\tothe{2}}$.
For the RIV and OA algorithms, we selected $\gamma = 0.0015$, $\hat{S}_{0} = \SI{10}{\micro\coulomb\tothe{2}}$ and $\mu_e = \SI{e-5}{\micro\coulomb\tothe{2}}$.
For adequate noise rejection, $C_\eta = 440 C$ for the RIV algorithm and $C_\eta = 380 I$ for the OA algorithm.
The sNDTFT window was fixed to $N_f = 24$ in all cases.

The selection of the parameters was based on the following heuristic: we selected a sufficiently small value for $N_f$ that is comparable to the period of the lowest expected resonance, and by trial and error, we found the values of $\gamma$ that provide an acceptable system response.
We increased $N_f$ to improve the algorithm response and partially mitigate noise.
In the case of the RIV and OA algorithms, we increased the value of $c_5$ from the initial value of \num{1} to improve the rejection of noise or decreased the value for a faster response.
Then, we selected $\mu_e$ to ensure that $\hat{S}_k$ does not become singular.
We always set $\hat{S}_{0}$ at a high value to avoid transients at the beginning of the simulation.

\begin{figure}
\centering
\includegraphics{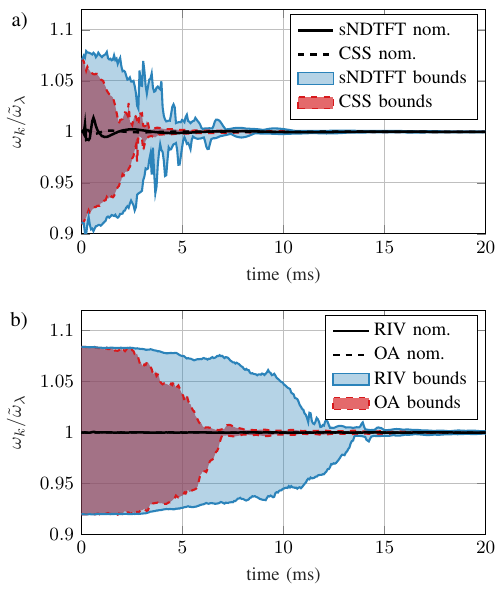}
\caption{Tracking of the resonance frequency for a set of 100 plants with 10\% parametric uncertainty. The region bounding the trajectory of the normalized frequency $\omega_k / \tilde{\omega}_\lambda$ for the set of plants is indicated by the shaded regions, along with the trajectory for the nominal plant. a) Simulation results using the cRPEM (the solid lines correspond to the LTI model with the sNDTFT algorithm, and the dashed lines correspond to the real plant in the CSS representation). b) Simulation results using the RIV and OA algorithms with the LTI model and the sNDTFT filter.}
\label{fig:piezoRes_direct}
\end{figure}

The bounds of the trajectories of $\omega_k / \tilde{\omega}_\lambda$ for the cRPEM algorithm are shown as blue solid lines in panel (a) of Fig.~\ref{fig:piezoRes_direct}.
Panel (b) of the same figure shows the bounds for the RIV and OA updates.
The cRPEM exhibits the fastest convergence, followed by the OA and the RIV updates.

To assess the effect of the sNDTFT, we repeated the simulations with the plant model transformed in the CSS representation, where $q_k$ is readily available.
Panel (a) in Fig.~\ref{fig:piezoRes_direct} compares the simulation results for the two system representations and the cRPEM algorithm, for which we observed the largest discrepancy.
The sNDTFT algorithm can be satisfactorily combined with the tracking algorithms.

\subsubsection{MIMO system of a gyroscope}
Furthermore, we simulated the tracking of the resonance frequency of a vibrating structure gyroscope \cite{Raman2009, Fei2015}.
The gyroscope contains two proof masses vibrating in a plane.
When the structure is rotated perpendicular to the plane, energy is transferred between the proof masses due to the Coriolis effect.
The vibratory gyroscope can be modeled as two 2\textsuperscript{nd}-order spring-mass-damper oscillators that are coupled by cross-damping and spring terms.
The cross-terms include the Coriolis effect and parasitic mechanical and electrical coupling.
The usual technique for acquiring the rotational speed is to excite one of the oscillators (which is referred to as the primary mode) at a constant amplitude while keeping the other oscillator (secondary mode) fixed.
As a result, the Coriolis effect acts on the secondary mode with a force that is proportional to the oscillating amplitude of the primary mode \cite{Raman2009}.
However, to measure the rotational speed accurately, the parasitic coupling must be either eliminated or identified and corrected.
Here, we propose an alternative approach to acquire the rotational speed acting on the gyroscope.
Both oscillators are excited at the same frequency, which coincides with the resonance of the primary mode.
The Coriolis effect alters the resonance frequency of the system, which in turn allows the rotational speed to be determined.

The input to our model is the control force $u_g = [u_p, u_s]^T$ that can be exerted on the primary and secondary modes.
The subscripts $p$ and $s$ indicate variables of the primary and secondary modes, respectively.
Given the displacement of the oscillators $x_g = [x_p, x_s]^T$, the system dynamics are described by
\begin{equation}\label{eq:gyro_ssmodel}
\begin{bmatrix} \dot{x}_g \\ \ddot{x}_g \end{bmatrix} =
\begin{bmatrix} 
0 & I \\
-K_g & -(D_g + \Omega) \\
\end{bmatrix} 
\begin{bmatrix} x_g \\ \dot{x}_g \end{bmatrix}
 + \begin{bmatrix} 0 \\ I \end{bmatrix} u_g
\end{equation}
where $D_g$ and $K_g$ are the damping and stiffness matrices of the system, respectively.
The Coriolis acceleration acts on the system through
\begin{equation}
\Omega = 
\begin{bmatrix*}
0 & -\omega_z \\
\omega_z & 0 \\
\end{bmatrix*}
\end{equation}
where $\omega_z$ is the rotational speed to be measured.
For our simulations, the model parameters were set to
\begin{IEEEeqnarray*}{c+c}
K_g =
\begin{bmatrix}
355.3 & 70.99\\
70.99 & 532.9 \\
\end{bmatrix} &
D_g =
\begin{bmatrix*}
0.01 & 0.002 \\
0.002 & 0.01 \\
\end{bmatrix*}
\end{IEEEeqnarray*}
as proposed in \cite{Fei2015}.
The model is normalized, and all of the units are dropped in the following.

The rotation $\omega_z$ does not induce an adequate resonance shift in the current gyroscope design.
The maximum shift of the primary mode can be achieved when $\dot{x}_s \approx x_p$.
The velocity of the secondary mode can be matched to the displacement of the primary mode by an LQR.
The LQR design for the CSS representation is described in Appendix~\ref{app:optimalLQR}.
To calculate the state feedback gain, $K$, we discretized (\ref{eq:gyro_ssmodel}) with $T_s = \SI{0.01}{\second}$ and converted the dicretized model into the CSS representation.
We solve the optimal control problem at the nominal working point, namely, $\omega_z = 0$.
We selected the state and input weights as
\begin{equation*}
Q_c = \begin{bmatrix} 
100 & 0 & 0 & -100 \\
0 & 0.001 & 0 & 0 \\
0 & 0 & 0.001 & 0 \\
-100 & 0 & 0 & 100
\end{bmatrix}
\end{equation*}
and $R_c = 10^{-4} I$.
The input to the gyroscope model is then synthesized as $s_k = K z_k + s_r$, with $z_k$ representing the complex envelope of the discretized states and $s_r$ representing a constant excitation.
We note that since $Q_c$ and $R_c$ have no imaginary part, $K$ also does not have an imaginary part.
For the resonance tracking algorithm, we consider (\ref{eq:gyro_ssmodel}) with the LQR feedback given the input $s_r$ and output $x_g$.

\begin{figure}
\centering
\includegraphics{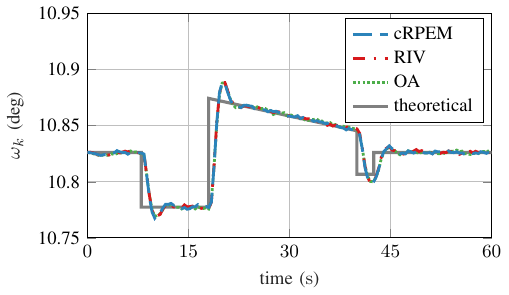}
\caption{Tracking of the resonance frequency of a vibratory gyroscope.}
\label{fig:gyro_sim}
\end{figure}

To assess the robustness of the different algorithms, we estimated $\mathcal{D}_c$ by solving (\ref{eq:optim_LyapReg}).
The OA algorithm was found to have the $\mathcal{D}_c$ with the largest area, with $d_m = 0.06$, $d_M = 1.02$ and $d_\phi = \ang{4.21}$.
The next largest area was identified for the RIV update, with $d_m = 0.71$, $d_M = 1.01$ and $d_\phi = \ang{0.94}$.
For the cRPEM, we found that $\mathcal{D}_c$ is demarcated by $d_m = 0.99$, $d_M = 1.02$ and $d_\phi = \ang{4.80}$.

Since the real system is not of the form (\ref{eq:lti_ch_state}-\ref{eq:lti_ch_output}), there is an offset in the estimated resonance that depends on $s_r$.
The offset for different values of $\omega_z$ can be found as the root of (\ref{eq:ode_f_closed}) from Theorem~\ref{th:closed_general}.
The input $s_r = [10, \; 18.6]^{T}$ eliminates the offset for all of the algorithms.

In the simulation of the gyroscope, we set the disturbance and measurement noise variances to $Q = \num{4e-4} I$ and $R = \num{1.6e-7} I$, respectively.
The estimation parameters were chosen following the heuristic procedure described in Section~\ref{sec:subsec_piezo}, resulting in $\gamma = 0.015$, $\hat{S}_0 = \num{e3}$ and $\mu_e = \num{8e-4}$ for the cRPEM and $\gamma = 0.0065$, $\hat{S}_0 = 50$ and $\mu_e = \num{e-6}$ for RIV and OA algorithms.
We selected $C_\eta = 20 C$ for RIV and $C_\eta = 20 I$ for OA.
In all cases, we set $N_f = 32$.
We simulated the response of the tracking algorithms assuming step and ramp changes in $\omega_z$.
The results are shown in Fig.~\ref{fig:gyro_sim}.
All of the algorithms yield similar results.
The performance is satisfactory, even for systems with rapidly changing parameters, confirming the effectiveness of the proposed scheme.

\section{Conclusion}\label{sec:conclusion}
In this work, we described a model-based resonance frequency tracking algorithm for linear systems.
We introduced a state transformation of linear systems into a complex-valued representation and lumped the resonance shift into a single variable.
This complex transformation allowed us to transfer well-established methods from the system identification framework and adjust the methods to analyze a general recursive algorithm for the current problem.
We described several special versions of the resonance tracking algorithm and examined their convergence.
We further discussed the implementation issues arising from the computation of the complex envelope of the signal in real time, and we validated our claims using numerical simulations.

Future work will include assessing alternative frequency updates to (\ref{eq:omega_update}) and modifying the presented algorithms to track multiple resonances by increasing the dimensions of $h$.
Moreover, the circularity of the complex noise and disturbance variables can be relaxed, extending the application of the tracking scheme to systems with widely linear complex random variables.

\section*{Acknowledgements}
The author would like to thank Dario Izzo (ESA, Advanced Concepts Team) for his valuable input and comments.

\bibliographystyle{unsrt}
\bibliography{ms}

\section*{Appendices}
\renewcommand{\setthesubsection}{\Alph{subsection}}
\begin{subappendices}
\subsection{Continuous-time CSS transformation}\label{app:continuousTimeCSS}
Consider the state-space representation of the continuous-time LTV system
\begin{IEEEeqnarray}{rCl}
\dot{x}(t) & = & \tilde{A}(t) x(t) + \tilde{B}(t) u(t) + \tilde{w}(t) \label{eq:lti_cont_state}\\
y(t) & = & \tilde{C}(t) x(t) + \tilde{D}(t) u(t) + \tilde{v}(t) \label{eq:lti_cont_output}
\end{IEEEeqnarray}
where $\dot{x}(t)$ is the time derivative of $x(t)$ with respect to $t$.
For convenience, we drop the notation $x(t)$ in favor of $x$ and use the subscript $i$ to indicate the $i$\textsuperscript{th} element of the corresponding vector.
We assume that $u_i = \beta_i e^{j (\theta +  \psi_i)}$ and $x_i = \alpha_i e^{j (\theta + \phi_i)}$.
We introduce the complex envelope variables with the components $z_i = \alpha_i  e^{j \phi_i}$ and $s_i = \beta_i e^{j \psi_i}$ and
compute the time derivative of the state variables:
\begin{IEEEeqnarray}{rCl}
\dot{z}_i & = & \dot{\alpha}_i e^{j \phi_i} + j \alpha_i \dot{\phi}_i e^{j \phi_i} \\
\dot{x}_i & = & \dot{\alpha}_i e^{j ( \theta + \phi_i)} + j \alpha_i ( \omega + \dot{\phi}_i) e^{j ( \theta + \phi_i)} = \dot{z}_i e^{j \omega t} + j \omega z_i  e^{j \omega t}
\end{IEEEeqnarray}
where $\omega = \dot{\theta}$.
Substitution into~(\ref{eq:lti_cont_state}-\ref{eq:lti_cont_output}) and elimination of the $e^{j \theta}$ terms results in the continuous-time CSS representation:
\begin{IEEEeqnarray}{rCl}
\dot{z}  & = & (\tilde{A} - j \omega I) z + \tilde{B} s + w \label{eq:css_cont_state} \\
q & = & \tilde{C} z + \tilde{D} s + v \label{eq:css_cont_output}
\end{IEEEeqnarray}
where $q$ is the complex envelope of $y$.
The noise is also transformed into its complex equivalent, similar to the discrete-time case.

The similarities between the transformed and original models discussed in the discrete-time case also apply to continuous-time models.
We note that only in the case of LTI systems, the zero-order hold discretization of the CSS-transformed system results in the discrete-time CSS representation.
For a constant $\omega$ in a time interval of length $T_s$, the matrix exponential of $A - j\omega I$ is
\begin{equation}
e^{(A - j\omega I) T_s} = e^{-j \omega T_s} e^{A T_s}
\end{equation}
since the matrices $A$ and $j\omega I$ commute.
The discretization of the complex matrix term $ -j\omega I$ results in the multiplication by $e^{-j \omega T_s}$ seen in the discrete-time CSS system.

\subsection{Optimal control for the CSS model}\label{app:optimalLQR}
Consider the optimal control problem with state update~(\ref{eq:linear_complex_state}) and the quadratic cost function:
\begin{equation}
J_c = E \left[ z_{N_c}^H Q_{c,N_c} z_{N_c} + \sum_{i = 0}^{N_c-1}   z_i^H Q_c z_i +  s_i^H R_c s_i \right]
\end{equation}
where $Q_{c,N_c}$, $Q_c$ and $R_c$ are real positive definite matrices of appropriate dimensions that penalize the terminal cost, state and control input, respectively.
The trajectory of the optimal cost can be computed by applying the dynamic programming algorithm \cite{Bertsekas1995} starting from the final cost:
\begin{equation}\label{eq:optCostN}
J^*_{c,N_c}(z_{N_c}) = z_{N_c}^H Q_{s,N_c} z_{N_c}
\end{equation}
where $J^*_{c,k}(z_k)$ denotes the optimal cost at time $k$ from $z_k$.
Similarly, at time $N_c-1$, the optimal cost is given as
\begin{IEEEeqnarray}{l}
J^*_{c,N_c-1}(z_{N_c-1}) = \min_{s_{N_c-1}} E \bigg[ z_{N_c-1}^H Q_c z_{N_c-1} +  s_{N_c-1}^H R_c s_{N_c-1} \IEEEnonumber \\
\quad + \> J^*_{c,N_c}\left( (\tilde{A}_{N_c-1} z_{N_c-1} \tilde{B}_{N_c-1} s_{N_c-1} + w_{N_c-1}) e^{-j \omega_k} \right) \bigg]. \IEEEeqnarraynumspace \label{eq:lqr_dpAlgNm1}
\end{IEEEeqnarray}
By setting the derivative of $J^*_{c,N_c-1}(z_{N_c-1})$ with respect to $s_{N_c-1}$ to zero, we recover the optimal input:
\begin{IEEEeqnarray}{rCl}
s^*_{N_c-1} & = & -K_{N_c-1} z_{N_c-1} \\
K_{N_c-1} & = & (R_c + \tilde{B}_{N_c-1}^H Q_{c,N_c} \tilde{B}_{N_c-1})^{-1} \tilde{B}_{N_c-1}^H Q_{c,N_c} \tilde{A}_{N_c-1}.
\end{IEEEeqnarray}
Substitution of the optimal input into ~(\ref{eq:lqr_dpAlgNm1}) results in
\begin{equation}
J^*_{c,N_c-1}(z_{N_c-1}) = z_{N_c-1}^H V_{N_c-1} z_{N_c-1} + E \left[ w_{N_c-1}^H Q_{c,N_c} w_{N_c-1} \right]  \label{eq:optCostNm1}
\end{equation}
where the optimal cost is quadratic with respect to the current state $z_{N_c-1}$.
The symmetric matrix $V_{N_c-1}$ is equal to
\begin{equation}
V_{N_c-1} = \tilde{A}_{N_c-1}^H Q_{c,N_c} \tilde{B}_{N_c-1} K_{N_c-1} + \tilde{A}_{N_c-1}^H Q_{c, N_c} \tilde{A}_{N_c-1}  + Q_c.
\end{equation}
The recursive application of the dynamic programming algorithm results in a quadratic representation of the optimal cost.
The weight matrix of the cost is given by the recursion
\begin{equation}
V_{k-1} = \tilde{A}_k^H V_k \tilde{B}_k (R_c + \tilde{B}_k^H V_k \tilde{B}_k)^{-1} \tilde{B}_k^H V_k \tilde{A}_k + \tilde{A}_k^H V_k \tilde{A}_k + Q_c  \label{eq:lrqCostMatrix}
\end{equation}
with the terminal value $V_{N_c} = Q_{c,N_c}$.
The optimal cost weighting matrix for the model~(\ref{eq:linear_uncertain_state}-\ref{eq:linear_uncertain_output}) follows the same recursion; for equal terminal costs, the trajectories of the optimal cost for the two models are identical.

\subsection{Duality between $\mathbb{CR}$ and real derivatives and norms}\label{app:complex_derivatives}
Here, we present the connection between the real-valued and $\mathbb{CR}$ derivatives and norms.
We use the Euclidean norm, which is defined by the inner product, and the matrix norm that is induced.
Given the vectors $h \in \mathbb{C}^n$ and $h_r = [ \Re h^T \; \Im h^T ]^T \in \mathbb{R}^{2n}$ and $h_c = [ h^T \; h^H ]^T \in \mathbb{C}^{2n}$, we can write \cite{KreutzDelgado2009, Dini2012}
\begin{equation}
\begin{bmatrix} h \\ \overline{h} \end{bmatrix} = 
\begin{bmatrix} I & \phantom{-}j I \\ I & -j I \end{bmatrix}
\begin{bmatrix} \Re h \\ \Im h \end{bmatrix}
\end{equation}
or $h_c = J h_r$.
Moreover, $h_r = \frac{1}{2} J^H h_c$.
Note that $\frac{1}{\sqrt{2}} J$ is unitary and that multiplication by a unitary matrix does not affect the norm.
Therefore,
\begin{equation}
\lVert h_r \rVert = \frac{1}{2} \lVert J^H h_c \rVert = \frac{\sqrt{2}}{2} \lVert h_c \rVert = \lVert h \rVert
\end{equation}
as expected.

We apply the same reasoning for the Jacobian matrices.
Given a function $f(h) \colon \mathbb{C}^n \mapsto \mathbb{C}^p$, the function can be written as $f_r(h_c) = f_r(h_r) = [ \Re f^T \; \Im f^T ]^T$ or as $f_c(h_c) = [ f^T \; f^H ]^T$.
The Jacobian matrices are related as \cite{KreutzDelgado2009}
\begin{IEEEeqnarray*}{c}
\nabla_r f_r = \nabla_c f_r J = \frac{1}{2} J^H \nabla_c f_c J \\
\nabla_r f_r = \begin{bmatrix}
\frac{\strut\partial \Re f}{\strut\partial \Re h} & \frac{\strut\partial \Re f}{\strut\partial \Im h} \\
\frac{\strut\partial \Im f}{\strut\partial \Re h} & \frac{\strut\partial \Im f}{\strut\partial \Im h}
\end{bmatrix}, \quad
\nabla_c f_r = \begin{bmatrix}
\frac{\strut\partial \Re f}{\strut\partial h} & \frac{\strut\partial \Re f}{\strut\partial \overline{h}} \\
\frac{\strut\partial \Im f}{\strut\partial h} & \frac{\strut\partial \Im f}{\strut\partial \overline{h}}
\end{bmatrix} \\
\nabla_c f_c = \begin{bmatrix}
\frac{\strut\partial f}{\strut\partial h} & \frac{\strut\partial f}{\strut\partial \overline{h}} \\
\frac{\strut\partial \overline{f}}{\strut\partial h} & \frac{\strut\partial \overline{f}}{\strut\partial \overline{h}}
\end{bmatrix}
\end{IEEEeqnarray*}
which means that $\lVert \nabla_r f_r \rVert = \sqrt{2} \lVert \nabla_c f_r \rVert = \lVert \nabla_c f_c \rVert$.
If $f$ is holomorphic, then $\lVert \nabla_r f_r \rVert = \lVert \frac{\partial f}{\partial h} \rVert$ since $\nabla_c f_c$ becomes block diagonal with blocks of equal norm.
Moreover, if $f(h) \colon \mathbb{C}^n \mapsto \mathbb{R}^p$, then for any given vector $z = [z_1^T \; z_2^T]^T$ with $\lVert z \rVert > 0$ and $z_1, z_2 \in \mathbb{C}^n$, $\nabla_r f z = \frac{\partial f}{\partial h} z_1 + \frac{\partial f}{\partial \overline{h}} z_2$ and
\begin{equation}
\lVert \nabla_r f z \rVert \leq \left \lVert \frac{\partial f}{\partial h} z_1 \right \rVert + \left \lVert \frac{\partial f}{\partial \overline h} z_2 \right \rVert \leq 2 \left \lVert \frac{\partial f}{\partial h} \right \rVert \lVert z \rVert
\end{equation}
since $\lVert z \rVert \geq \lVert z_1 \rVert$, $\lVert z \rVert \geq \lVert z_2 \rVert$ and $\lVert \frac{\partial f}{\partial h} \rVert = \left \lVert \overline{\frac{\partial f}{\partial h}} \right \rVert = \lVert \frac{\partial f}{\partial \overline{h}} \rVert $.

\subsection{Proof of Theorem~\ref{th:openloop}}\label{app:proofThopen}

The convergence properties of the algorithm (\ref{eq:rpem_est_output}-\ref{eq:rpem_state_der}) follow from \cite[theorem 2]{Ljung1981} after verifying the necessary regularity conditions (labeled Conditions S1, M1 and A1 in \cite{Ljung1981} and described in the following).
Since the excitation frequency can be arbitrary, it is convenient to incorporate $\omega_k$ in the input, similar to (\ref{eq:lti_complex_envelope}).
We apply a change of variables to the estimator model:
\begin{IEEEeqnarray}{rCl}
\hat{x}_{k+1} & = & \hat{h}_{k+1} [ (A - L C) \hat{x}_k  + (B - L D) u_k + L y_k ] \IEEEeqnarraynumspace \\
\xi_{k+1} & = & \frac{\hat{x}_{k+1}}{\hat{h}_{k+1}} + \hat{h}_{k+1} (A - L C) \xi_{k}
\end{IEEEeqnarray}
with $\hat{x}_k = \hat{z}_k e^{j \theta_k}$, $u_k = s_k e^{j \theta_k}$, $y_k = q_k e^{j \theta_k}$ and $\xi_k = \eta_k e^{j \theta_k}$.
The $\tilde{h}_k$ and $\tilde{S}_k$ updates (\ref{eq:rpem_hest}-\ref{eq:rpem_hess_approx}) and the definition of $J_{est,a}$ are not affected by the change of variables since $e^{j \theta_k}$ is counteracted by its conjugate.
The resulting $\mathcal{M}$ is an LTI model set for a fixed $h$.
All of the models in $\mathcal{M}$ are twice differentiable with respect to $h$, and their eigenvalues lie strictly inside the unit circle $\forall h \in \mathcal{D}_\mathcal{M}$, satisfying Condition~M1.

Condition~S1 requires the data generation of the real system to be exponentially stable in the sense that, for each $k, l \colon k \geq l$, $\exists \, y_k^0, u_k^0$ independent of $y_l, u_l$ such that $E (\lVert y_k - y_k^0 \rVert + \lVert u_k - u_k^0 \rVert)^8 < c_3 \lambda_0^{k-l}, c_3 < \infty, \lambda_0 < 1$.
For the real system, we have used the equivalent representation (\ref{eq:lti_complex_envelope}) with $x_k = z_k e^{j \theta_k}$.
Starting from the initial state $x_l$, we have
\begin{equation}
x_k = \Phi(k,l) x_l + \sum_{i=l}^{k-1} \Phi(k,i+1) (\tilde{B} u_i + w_i)
\end{equation}
and $y_k = \tilde{C} x_k + \tilde{D} u_k + v_k$.
Then, $x_k^0$ can be estimated as the second part of the previous equation \cite{Ljung1978}, which does not depend on $x_l$.
The input is independent of past values, $u_k^0 = u_k$.
Therefore,
\begin{equation}
E \big [ \lVert y_k - y_k^0 \rVert^8 \big ] = E \big [ \lVert C (x_k - x_k^0) \rVert^8 \big ] = \lVert C \Phi(k,l) x_l \rVert^8 \leq \lVert C \rVert^8 \lVert \Phi(k,l) \rVert^8 \lVert z_l \rVert^8
\end{equation}
satisfies the definition as a result of Lemma~\ref{lem:stb}.

For fixed $h$, the sequence of $\rho_k e^{j \theta_k}$ is bounded since the sequences $y_k$ and $u_k$ are bounded and the system is stable.
The covariance of $\rho_k$ is also bounded.
Therefore, the limit (\ref{eq:estObjFunInf}) is well defined, which satisfies Condition~A1.
Finally, the necessary requirement for $\gamma_k$ in \cite{Ljung1981} is identical to Assumption~\ref{assum:gamma_gain}.

With the conditions verified, according to \cite[Theorem 2]{Ljung1981}, the convergence of the algorithm is subject to the asymptotic stability of the associated ODE, and the RPEM converges with probability \num{1} to a local minimum of $J_{est,a}$ or to the boundary of $\mathcal{D}_\mathcal{M}$ as $k \to \infty$.

\subsection{Sliding nonuniform discrete-time Fourier transform}\label{app:sDFT_der}
The conversion of the measured real-valued signal into the complex envelope representation can be accomplished by applying the sliding nonuniform discrete-time Fourier transform (sNDTFT). 
Nonuniformity refers to the fact that the instantaneous frequency of the signal may not be constant (although it is known).
The complex envelope $q_k$ can be approximated by an sNDTFT of length $N_f$, $Y_k$, of the measured signal $y_k$ as follows:
\begin{equation}
Y_k = y_{k-N_f+1} + \sum_{i=2}^{N_f} y_{k-N_f+i} e^{-j \sum_{l = 2}^i \omega_{k-N_f+l}} = \sum_{i=0}^{N_f-1} y_{k-i} e^{-j \delta \theta_k(i,N_f-1)} \label{eq:ndft_basic}
\end{equation}
where $\delta \theta_k(i,l) = \theta_{k-i+1} - \theta_{k-l+1}$ is the phase difference between the samples.
Following \cite{Sysel2012,Jacobsen2003}, we derive a recursive method to compute~(\ref{eq:ndft_basic}).
First, we multiply both sides of~(\ref{eq:ndft_basic}) by $e^{j \delta \theta_k(0,N_f-1)}$:
\begin{equation}
Y_k = e^{-j \delta \theta_k(0,N_f-1)} \sum_{i=0}^{N_f-1} y_{k-i} e^{j \delta \theta_k(0,i)} = e^{-j \delta \theta_k(0,N_f-1)} \tilde{Y}_k.
\end{equation}
The second term can be computed recursively as
\begin{equation}
 \tilde{Y}_k =  \tilde{Y}_{k-1} e^{j \omega_k} - y_{k-N_f}  e^{j \delta \theta_k(0,N_f)} + y_k.
\end{equation}
The value of $\delta \theta_k(0,N_f)$ can be updated at each time step, as in~(\ref{eq:sdft_binRemAngle}).
Phase and magnitude corrections must be applied to $Y_k$ to recover $q_k$.
The previous calculation of the sNDTFT assumes zero phase at the start of the computation window, so we must offset the calculation by $\theta_{k-N_f+1}$ to be consistent when comparing the phase shifts to $\theta_k$:
\begin{equation}
q_k = 2 \frac{Y_k}{N_f} e^{j \theta_{k-N_f+1}} = 2 \frac{ \tilde{Y}_k}{N_f} e^{j \theta_k}.
\end{equation}

\end{subappendices}

\end{document}